\documentclass[11pt]{amsart}
\usepackage{enumerate, amsmath, amsfonts, amssymb, amsthm, mathtools, thmtools, wasysym, graphics, graphicx, xcolor, frcursive,xparse,comment,bbm}
\usepackage[colorlinks=true, pdfstartview=FitV, linkcolor=blue, citecolor=blue, urlcolor=blue]{hyperref}
\usepackage{url, hypcap}
\hypersetup{colorlinks=true, citecolor=darkblue, linkcolor=darkblue}
\definecolor{darkblue}{rgb}{0.0,0,0.7} % darkblue color
\newcommand{\AG}{\mathcal{A}_G}
\usepackage[colorinlistoftodos]{todonotes}

\usepackage{wasysym, MnSymbol, hyperref}
\usepackage{ytableau}
\usepackage{tikz-cd} 
\usepackage{graphics}

\newcommand{\Tab}{\mathrm{Tab}}
\newcommand{\wt}{\mathsf{wt}}
\newcommand{\sh}{\mathsf{sh}}
\newcommand{\stand}{\mathsf{stand}}
\usepackage[margin= 1 in]{geometry}

\newtheorem{theorem}{Theorem}[section]
\newtheorem{lemma}[theorem]{Lemma}
\newtheorem{remark}[theorem]{Remark}
\newtheorem{example}[theorem]{Example}
\newtheorem{proposition}[theorem]{Proposition}
\newtheorem{corollary}[theorem]{Corollary}
\theoremstyle{definition}
\newtheorem{definition}[theorem]{Definition}
\newtheorem{problem}[theorem]{Problem}
\numberwithin{equation}{section}

\usepackage{subcaption}
\definecolor{darkred}{rgb}{0.7,0,0} % darkred color
\newcommand{\defn}[1]{{\color{darkred}\emph{#1}}} % emphasis of a definition

\author[J.~Pappe]{Joseph Pappe}
\address[J. Pappe]{Department of Mathematics, UC Davis, One Shields Ave., Davis, CA 95616-8633, U.S.A.}
\email{jhpappe@ucdavis.edu}

\author[D.~Paul]{Digjoy Paul}
\address[D. Paul]{Chennai Mathematical Institute, Chennai, India.}
\email{digjoypaul@gmail.com}

\author[A.~Schilling]{Anne Schilling}
\address[A. Schilling]{Department of Mathematics, UC Davis, One Shields Ave., Davis, CA 95616-8633, U.S.A.}
\email{anne@math.ucdavis.edu}

\date{\today}
\title{The Burge correspondence and crystal graphs}

\begin{document}

\begin{abstract}
The Burge correspondence yields a bijection between simple labelled graphs and semistandard Young tableaux of threshold shape.
We characterize the simple graphs of hook shape by peak and valley conditions on Burge arrays. This is the first step towards an
analogue of Schensted's result for the RSK insertion which states that the length of the longest increasing subword of a word is the length 
of the largest row of the tableau under the RSK correspondence. Furthermore, we give a crystal structure on simple graphs of hook shape. 
The extremal vectors in this crystal are precisely the simple graphs whose degree sequence are threshold and hook-shaped.
\end{abstract}

\maketitle

%%%%%%%%%%%%%%%%%%%%%%%%%%%%%%%%%%%%%%%%%%%%%%%
\section{Introduction}

The celebrated \defn{Robinson--Schensted (RS) correspondence}~\cite{Robinson,Schensted.1961} gives a bijection between
words $w$ in the alphabet $\{1,2,\ldots,n\}$ of length $k$ and a pair of tableaux of the same shape $\lambda$, a partition of $k$ with at 
most $n$ parts, where the first tableau is a semistandard Young tableau in the same alphabet and the second tableau is a standard tableau. 
Schensted~\cite{Schensted.1961} proved that $\lambda_1$ (the biggest part of the partition $\lambda$) is the length of the longest increasing subword of 
$w$. Knuth's generalization of the RS correspondence~\cite{Knuth}, known as the RSK correspondence, provides a bijective proof of the \defn{Cauchy identity} 
in symmetric function theory
\[
	\sum_{\lambda} s_\lambda(x) s_\lambda(y) = \prod_{i,j \geqslant 1} \frac{1}{1-x_i y_j},
\]
where the sum is over all partitions $\lambda$ and $s_\lambda(x)$ is the Schur function in the variables $x_1,x_2,\ldots$ indexed by the partition $\lambda$.

In~\cite{Burge.1974}, W. Burge gives four variants of the RSK correspondence. In this paper, we focus on the correspondence 
in~\cite[Section 4]{Burge.1974}, which gives a bijection between simple labelled graphs (graphs without loops or multiple edges) and semistandard Young 
tableaux of threshold shape. A partition $\lambda=(\lambda_1,\lambda_2,\ldots, \lambda_n)$ is called \defn{threshold} if $\lambda^t_{i}=\lambda_i +1$ for 
all $1 \leqslant i \leqslant d(\lambda)$, where $\lambda^t_{i}$ is the length of $i$-th column of the Young diagram of $\lambda$ and $d(\lambda)$
is the maximal $d$ such that $(d,d) \in \lambda$. We call this bijection the \defn{Burge correspondence}.
The Burge correspondence gives a bijective proof of the \defn{Littlewood identity}~\cite[Exer. I.5.9(a) and I.8.6(c)]{MR3443860}
\begin{equation}\label{littlewood}
	1+\sum_{\lambda} s_{\lambda}(x_1,x_2,\ldots) = \prod_{i<j} \left(1+x_ix_j\right),
\end{equation}
where the sum runs over all threshold partitions. For more details about the representation theoretic significance of the Burge correspondence, 
see~\cite{PPS.2019}.

In this paper, we characterize the graphs whose shapes under the Burge correspondence are hook shapes in terms of \defn{peak} and 
\defn{valley} conditions. This is the first step towards an analogue for the Burge correspondence of Schensted's result for the RS correspondence,
namely that increasing sequences under the RS correspondence give tableaux of single row shape.

\medskip

Threshold partitions also play an important role in graph theory.
Given a simple graph $G=([n],E)$ with vertex set $[n]=\{1,2,\ldots,n\}$ and edge set $E$, the degree $d_i$ of a vertex $i$ is the number 
of neighbors of $i$. The \defn{degree sequence} of $G$ is the tuple $d_G=(d_1,d_2,\ldots,d_n)$. 
One of the central questions in graph theory is the characterization of all sequences that appear as degree sequences of a simple graph
(see for example~\cite{ErdosGallai.1960,RuchGutman.1979,PeledSrinvasan.1989,MP.1995,MerrisRoby.2005}).
The Erd\"os--Gallai theorem~\cite{ErdosGallai.1960} gives a characterization of graphic partitions, that is, partitions that are degree sequences of a simple graph. 

Define $X^G:=x_1^{d_1}x_2^{d_2}\cdots$, the monomial associated with a simple graph $G$ with degree sequence $d_G=(d_1,d_2,\ldots)$. 
It is known that the generating function $\sum_{G} X^G=\prod_{i<j} \left(1+x_ix_j\right)$ is a symmetric function in $x_1,x_2,\ldots$, where the sum runs 
over all simple graphs $G$. Note that the right hand side is the right hand side in Littlewood's identity~\eqref{littlewood}. Hence, from~\eqref{littlewood} 
and the Burge correspondence, it follows that a sequence $d=(d_1 \geqslant d_2 \geqslant \cdots \geqslant d_n)$ is 
a degree sequence of a simple graph $G$ if and only if $d \leqslant \lambda$ for some threshold partition $\lambda$, where $\leqslant$ is the 
dominance order for partitions. 
In \cite{MR1437004}, Gasharov proved that this condition is equivalent to the Erd\"os--Gallai theorem. Hence threshold partitions 
are the maximal graphic partitions.

The \defn{degree partition} of $G$ is the partition $\tilde{d}_G$ obtained by rearranging $d_G$ in weakly decreasing manner. 
A graph $G$ is called \defn{threshold} if the associated degree partition is threshold. For alternative definitions and characterizations of threshold graphs, 
see~\cite[Chapter 3]{MP.1995}.
Threshold graphs can be interpreted in an extremal sense. First, the threshold partitions (the degree partitions of threshold graphs) are maximal 
among graphic partitions. Second, the threshold partitions of $n$ are the extreme points of the degree partition polytope, the convex hull of all degree partitions 
of simple graphs on $n$, see \cite{MR2223521} and references therein. 

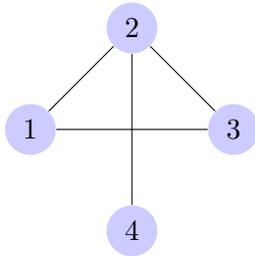
\begin{figure}[h]
\begin{center}
\begin{tikzpicture}  
  [scale=.9,auto=center, every node/.style={circle,fill=blue!20}]
    
  \node (a1) at (0,0) {1};  
  \node (a2) at (3,0)  {3};   
  \node (a3) at (1.5,1.5)  {2};  
  \node (a4) at (1.5,-1.5) {4};  
  
  \draw (a1) -- (a2); % these are the straight lines from one vertex to another  
  \draw (a1) -- (a3);  
  \draw (a2) -- (a3);  
  \draw (a3) -- (a4);  
\end{tikzpicture}
\caption{A threshold graph
\label{figure.threshold}} 
\end{center} 
\end{figure}

\begin{example}
The simple graph in Figure~\ref{figure.threshold} is threshold as its degree partition $(3,2,2,1)$ is threshold. 
\end{example}

The condition $\lambda^t_{i}=\lambda_i +1$ for threshold partitions is in fact Merris' and Roby's reformulation~\cite{MerrisRoby.2005} of
conditions stated by Ruch and Gutman~\cite{RuchGutman.1979}. Klivans and Reiner~\cite{KlivansReiner.2008} give some generalizations
of these concepts to hypergraphs and related them to plethysm.

In this paper, we characterize when the shape $\lambda_G$ of a simple graph $G$ is a hook shape. This is done by analyzing the Burge 
correspondence and imposing peak and valley conditions on the Burge array corresponding to $G$. In addition, we impose a crystal structure
on simple graphs of hook shape. Crystal graphs are combinatorial skeletons of Lie algebra representations (see for example~\cite{BumpSchilling.2017,
HongKang.2002}). 
The extremal vectors in this crystal are precisely the simple graphs whose degree sequence is threshold.

The paper is organized as follows. In Section~\ref{section.burge}, we review the Burge correspondence and prove some results
regarding standardization. In Section~\ref{section.shape}, we provide the characterization of simple graphs of hook shape. 
Finally in Section~\ref{section.crystal}, we give the crystal structure on simple graphs of hook shape.

%%%%%%%%%%%%%%%%%%%%%%%%%%%%%%%%%%%%%%%%%%%%%%%%%%
\subsection*{Acknowledgements}
We thank Amritanshu Prasad for fruitful discussions.

AS was partially supported by NSF grants DMS--1760329 and DMS--2053350.

%%%%%%%%%%%%%%%%%%%%%%%%%%%%%%%%%%%%%%%%%%%%%%%%%%
\section{The Burge correspondence}
\label{section.burge}
%%%%%%%%%%%%%%%%%%%%%%%%%%%%%%%%%%%%%%%%%%%%%%%%%%

In this section, we define the Burge correspondence~\cite{Burge.1974}. We review some preliminaries in Section~\ref{section.prelim}.
We remind the reader of Schensted's result on longest increasing subwords of words in Section~\ref{section.longest increasing}
before introducing the Burge correspondence in Section~\ref{section.sub Burge}. In Section~\ref{section.standardization}, we show that
the Burge correspondence intertwines with standardization.

%%%%%%%%%%%%%%%%%%%%%%%%%%%%%%%%%%%%%%%%%%%%%%%%%%
\subsection{Preliminaries}
\label{section.prelim}

A \defn{partition} $\lambda$ of a nonnegative integer $n$, denoted by $\lambda \vdash n$, is a weakly decreasing sequence $\lambda=(\lambda_1,\lambda_2,\ldots,\lambda_\ell)$ of positive
integers $\lambda_i$ such that $\sum_{i=1}^\ell \lambda_i =n$. The \defn{length} of $\lambda$ is $\ell$. The \defn{Young diagram} $Y(\lambda)$ of 
$\lambda$ is a left-justified array of boxes with $\lambda_i$ boxes in row $i$ from the top. (This is also known as the English convention for Young 
diagrams of partitions). A partition $\lambda$ is a \defn{hook} if $Y(\lambda)$ does not contain any $2\times 2$ squares.

\begin{definition}
Let $\lambda$ be a partition. A \defn{semistandard Young tableau} of shape $\lambda$ in the alphabet $[n]:=\{1,2,\ldots,n\}$ is a filling of the
Young diagram of $\lambda$ with letters in $[n]$ such that the numbers weakly increase along rows and strictly increase along columns.
We denote by $\Tab_n(\lambda)$ the set of all semistandard Young tableaux of shape $\lambda$ in the alphabet $[n]$.
\end{definition}

Let $T$ be a semistandard Young tableau. The \defn{shape} of $T$ is denoted $\sh(T)$.
The \defn{weight} of a semistandard Young tableau $T$, denoted $\wt(T)$, is the integer vector $(\mu_1,\ldots,\mu_n)$, where $\mu_i$ is the number 
of times the number $i$ occurs. The subset of $\Tab_n(\lambda)$ consisting of all semistandard Young tableaux of weight $\mu$ is denoted by 
$\Tab(\lambda,\mu)$.

Given an integer vector $\mu=(\mu_1,\dotsc,\mu_n)$, let $x^\mu$ denote the monomial $x_1^{\mu_1} x_2^{\mu_2} \cdots x_n^{\mu_n}$ in the $n$ variables
$x_1,x_2,\ldots,x_n$.

\begin{definition}
  For each integer partition $\lambda$, the \defn{Schur polynomial} in $n$ variables corresponding to $\lambda$ is defined as
  \begin{displaymath}
    s_\lambda(x_1,\dotsc,x_n) = \sum_{T \in \Tab_n(\lambda)} x^{\wt(T)}.
  \end{displaymath}
\end{definition}

%%%%%%%%%%%%%%%%%%%%%%%%%%%%%%%%%%%%%%%%%%%%%%%%%%
\subsection{Schensted algorithm and longest increasing subwords}
\label{section.longest increasing}

The Burge correspondence (as well as the celebrated Robinson--Schensted--Knuth (RSK) correspondence) uses the \defn{Schensted row insertion} algorithm. 
Given a semistandard Young tableau $T$ in the alphabet $[n]$, a letter $i \in [n]$ can be inserted into $T$ in the following way: if $i$ is larger than 
or equal to all the entries of the first row of $T$, a new box containing $i$ is added at the end of first row and the process stops. Otherwise, $i$ replaces 
the smallest leftmost number $j$ of the first row such that $j>i$. Then $j$ is inserted in the second row of $T$ in the same way and so on. The procedure stops 
when a new box is added to $T$ at the end of a row. The result is denoted $T\leftarrow i$.
The shape of $T\leftarrow i$ contains one new box compared to the shape of $T$. 

Given a tableau $T$ and a letter $x$ that lies at the end of some row of $T$, the \defn{Schensted reverse bumping} algorithm generates a pair of a tableau 
$T'$ and a letter $y$ in the following way:
Let $t$ be the row index of $x$ and $x_1$ be the rightmost entry of row $t-1$ such that $x_1<x$. Replace  $x_1$ by $x$ in $T$ and output $x_1$. 
Repeat the process for $x_1$ and continue until an element of the first row say $y$ is obtained as output. The resulting tableau is $T'$.
We shall denote the pair $(T',y)$ by $T\rightarrow x$.

Given a word $w=w_1w_2\ldots w_k$ in the alphabet $\{1,2,\ldots, n\}$, the \defn{Schensted insertion tableau} is defined as 
$P(w):=\emptyset \leftarrow w_1 \leftarrow w_2 \leftarrow \cdots \leftarrow w_k$. The shape of the semistandard Young tableau $P(w)$,
denoted $\lambda(w)=(\lambda_1,\lambda_2,\ldots)$, is called the \defn{shape of the word} $w$. Schensted~\cite{Schensted.1961} proved that 
$\lambda_1$ is the length of the longest increasing subword of $w$. In particular, the shape $\lambda(w)$ is a single row if $w$ is weakly increasing.
Greene \cite{MR354395} extended the result of Schensted by interpreting the rest 
of the shape of $\lambda$. For a poset theoretic viewpoint of the map $w \mapsto \lambda(w)$ and various applications including in the context of 
flag varieties, see Britz and Fomin \cite{MR1814900} and references therein. In the same spirit, we ask very similar questions about the Burge 
correspondence introduced in the next subsection.

%%%%%%%%%%%%%%%%%%%%%%%%%%%%%%%%%%%%%%%%%%%%%%%%%%
\subsection{The Burge correspondence}
\label{section.sub Burge}

We begin by recalling the definition of a Burge array from~\cite[Section 4]{Burge.1974}.

\begin{definition}
\label{burge array}
Given a simple graph $G=([n],E)$ with $|E|=r$, define a two line array known as 
the \defn{Burge array}
\[
\mathcal{A}_G = 
\begin{bmatrix}
    a_1       & a_2 & \dots & a_{r} \\
    b_1       & b_2 & \dots & b_r
\end{bmatrix}
\]
satisfying:
\begin{enumerate}
\item Each pair $ (a_k,b_k) $ is an edge of $G$ and $a_k >b_k$ for each $1\leqslant k \leqslant r$.
\item The top line is weakly increasing, that is, $a_k \leqslant a_{k+1}$ for all $1\leqslant k < r$.
\item If $a_k=a_{k+1}$ for some $1 \leqslant k < r$, then $b_k>b_{k+1}$.
\end{enumerate}
\end{definition}

Notice that $G$ is completely determined by the associated Burge array $\mathcal{A}_{G}$ assuming that $[n]$ is known. Note that singletons
in the simple graph do not appear in the Burge array $\AG$.

Given a threshold partition $\lambda$, by definition, the Young diagram $Y(\lambda)$ of $\lambda$  is divided into two symmetric pieces. The bottom 
piece consists of all boxes that lie strictly below the diagonal and the top piece consists of the rest. Each position in the top (bottom) piece of the Young 
diagram of $Y(\lambda)$ corresponds to a unique position, called the \defn{opposite position}, in the bottom (top) piece of $Y(\lambda)$. The opposite 
position $\mathsf{op}(s,t)$ of $(s,t)$ is defined to be $(t+1,s)$ if $s\leqslant t$ and $(t,s-1)$ otherwise.

Having defined all the necessary tools, we now state the main algorithm of the \defn{Burge correspondence}.  Starting with the empty tableau $T_0$, we shall 
insert all the edges of $G$, as ordered in $\AG$, into $T_0$. Let $T_k$ be the tableau obtained by inserting the edge $(a_k,b_k)$ into $T_{k-1}$ in the 
following way: 
\begin{enumerate}
\item
First insert $b_k$ into $T_{k-1}$ using the Schensted insertion algorithm. This adds a new cell to the shape, say in position $(s_k,t_k)$.
\item
Place the entry $a_k$ in the cell $\mathsf{op}(s_k,t_k)$. Observe that each addition of an edge transforms a tableau of threshold shape to another 
tableau of threshold shape. 
\end{enumerate}
Finally, the tableau $T_G:=T_r$ is the threshold tableau associated with the graph $G$ under the Burge correspondence.

Burge~\cite{Burge.1974} proved that this tableau is semistandard. Given such a tableau $T$, we recover the Burge array in the following way: 
Let $a_r$ be the largest entry of $T$ with largest column index. Remove $a_r$ from $T$. Let $z_r$ be the value at the opposite position of the cell
containing $a_r$. Let $b_r=y_r$ where $T\rightarrow z_r=(T_{r-1},y_r)$. Repeat the process for $T_{r-1}$ and continue until the empty tableau is 
obtained in the output. 

\begin{example}
\label{eg nonhook}
The Burge array $\mathcal{A}_G$ for the graph in Figure \ref{figure.threshold} is
 $ 
\begin{bmatrix}
   2 & 3 & 3 & 4 \\
   1 & 2 & 1 & 2
\end{bmatrix}
$. The associated tableau $T_G$ of threshold shape $(3,2,2,1)$ is obtained by inserting the edges of $G$ as ordered in $ \mathcal{A}_G$
as depicted in Figure~\ref{figure.Burge insertion}.

\begin{figure}[t]
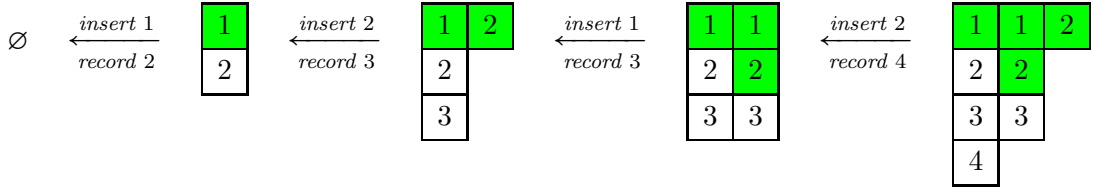

$\emptyset$
\quad
$\xleftarrow[\text{record } 2]{\text{insert }1}$
\quad
\begin{ytableau}
*(green) 1\\
2
\end{ytableau}
\quad
$\xleftarrow[\text{record } 3]{\text{insert }2}$
\quad
\begin{ytableau}
*(green) 1 & *(green) 2\\
2\\
3
\end{ytableau}
\quad
$\xleftarrow[\text{record } 3]{\text{insert }1}$
\quad
\begin{ytableau}
*(green) 1 & *(green) 1\\
2 & *(green) 2\\
3 & 3
\end{ytableau}
\quad
$\xleftarrow[\text{record } 4]{\text{insert }2}$
\quad
\begin{ytableau}
*(green) 1 & *(green) 1 & *(green) 2\\
2 & *(green) 2\\
3 & 3\\
4
\end{ytableau}
\caption{Burge insertion associated to Example~\ref{eg nonhook}.
\label{figure.Burge insertion}}
\end{figure}
\end{example}  

In the same spirit as the shape of a word under the RSK correspondence, we can define the shape of a graph under the Burge correspondence.

\begin{definition}
The \defn{shape of a simple graph} $G$ is the partition $\lambda_G:=\sh(T_G)$, where $T_G$ is the tableau associated with $G$ by the 
Burge correspondence.
\end{definition}

It can be observed that the shape of a threshold graph $G$ is its degree partition $\tilde{d}_G$. 
Namely, let $G$ be a threshold graph with degree sequence $d_G$. If $T$ is the semistandard Young tableau of threshold shape $\lambda$
and weight $d_G$, then $d_G$ is less than or equal to $\lambda$ in dominance order. Since $d_G$ is a threshold sequence by assumption,
we know that the only partition that dominates a threshold sequence is the corresponding partition. Hence $\lambda= \tilde{d}_G$.

We call a simple graph $G$ a \defn{hook-graph} if the associated tableau $T_G$ has hook shape. Given the nature of the Burge algorithm,
determining when a $T_G$ has hook shape is analogous to asking when a tableau under the RSK algorithm has single row shape.

\begin{problem}
\label{problem.shape}
What is the shape of a simple graph?
\end{problem}

In the next section, we characterize all hook-graphs.

%%%%%%%%%%%%%%%%%%%%%%%%%%%%%%%%%%%%%%%%%%%%%%%%%%%%%%%%%%%%%%%%%
\subsection{Standardization}
\label{section.standardization}

Both $\Tab_m(\lambda)$ and words over $[m]$ with length $n$ have the notion of \defn{standardization}. Standardization intertwines
with RSK in the sense that they form a commuting diagram. We show that an analogous result holds true for the Burge correspondence.

First, we review the standardization map for semistandard Young tableaux. Let $\lambda \vdash n$ and let $C = \{c_{1} < \cdots < c_{n}\}$ be a subset of 
$\mathbb{N}$. The standardization of $T \in \Tab(\lambda, \mu)$ with respect to the alphabet $C$, denoted by $\stand_{C}(T)$, is the map replacing all 
the $1$'s in $T$ from left to right with the numbers $c_{1}$ through $c_{\mu_{1}}$, replacing all the $2$'s from left to right with the numbers 
$c_{\mu_{1}+1}$ through $c_{\mu_{1}+\mu_{2}}$, etc.

The standardization map on words over the alphabet $[m]$ can be defined similarly. Let $\omega$ be a word using the alphabet $[m]$ with length 
$n$ and let $\mu$ denote its content. In other words, let $\mu = (\mu_{1}, \mu_{2}, \ldots, \mu_{m})$ be an integer vector, where $\mu_{i}$ denotes 
the number of $i$'s in $\omega$. The standardization of $\omega$ with respect to $C$, which we also denote by $\stand_{C}(\omega)$, is defined 
by replacing the $1$'s in $\omega$ from left to right with the numbers $c_{1}$ through $c_{\mu_{1}}$, replacing the $2$'s in $\omega$ from left to 
right with the numbers $c_{\mu_{1}+1}$ through $c_{\mu_{1}+\mu_{2}}$, etc. The following result formalizes the relationship between standardization 
and RSK, see for example~\cite[Lemma 7.11.6]{Stanley.EC2}.

\begin{proposition}
\label{prop.standardization}
Let $\omega$ be a word in the alphabet $[m]$ with length $n$. Then $\stand_{C}(P(\omega)) = P(\stand_{C}(\omega))$.
\end{proposition}

Analogously, we define the standardization of a Burge array $\AG$. Given a Burge array $\AG$ with $r$ columns and a subset 
$C= \{c_{1} < \cdots < c_{2r} \}$ of $\mathbb{N}$, define the standardization of $\AG$ with respect to $C$, denoted by $\overline{\stand}_{C}(\AG)$, 
to be the map that replaces the $1$'s in $\AG$ from left to right with the numbers $c_{1}$ through $c_{d_{1}}$, replaces the $2$'s from left to right with 
the numbers $c_{d_{1}+1}$ through $c_{d_{1}+d_{2}}$, etc. where $d_{G} = (d_{1}, d_{2}, \ldots, d_{n})$ is the degree sequence of $G$.

For $T$ a semistandard Young tableau, the \defn{reading word} of $T$ denoted by $R(T)$ is obtained by reading the entries within a row from left to 
right starting with the bottommost row.

\begin{proposition}
Let $\AG$ be a Burge array and $\mathcal{B}_{G} = \overline{\stand}_{C}(\AG)$. Let $T_{G}$ (resp. $S_{G}$) be
the tableau associated to $\AG$ (resp. $\mathcal{B}_{G}$) under the Burge correspondence. Then $S_{G} = \stand_{C}(T_{G})$.
\end{proposition}

\begin{proof}
We proceed by induction on the number of columns of 
$\AG =
\begin{bmatrix}
    a_1       & a_2 & \dots & a_{r} \\
    b_1       & b_2 & \dots & b_r
\end{bmatrix}$. 
Note that the base case of $r = 0$ is trivial. Let $\mathcal{A}_{r-1}$ (resp. $\mathcal{B}_{r-1}$) denote the Burge array formed by the first $r-1$ columns 
of $\AG$ (resp. $\mathcal{B}_{G}$). We have $\mathcal{B}_{r-1} = \overline{\stand}_{C-\{c_{i_{r}}, c_{2r}\}}(\mathcal{A}_{r-1})$ where $[c_{2r}, c_{i_{r}}]^{T}$ 
is the last column of $\mathcal{B}_{G}$. From our inductive hypothesis, $S_{r-1} = \stand_{C-\{c_{i_{r}}, c_{2r}\}}(T_{r-1})$. As the reading word $R(S_{r-1})$ 
is the standardization of $R(T_{r-1})$, we have $(R(S_{r-1}), c_{i_{r}})$ is the standardization of $(R(T_{r-1}), b_{r})$. By 
Proposition~\ref{prop.standardization}, $S_{r-1} \leftarrow c_{i_{r}}= \stand_{C-\{c_{2r}\}}(T \leftarrow b_{r})$. Thus the $a_{r}$ and $c_{2r}$ must be 
placed in the same position of their respective tableau. From the inverse of the Burge correspondence, $a_{r}$ is the largest value in $T_{G}$ and lies 
in a column to the right of any equivalent letters. Therefore, $a_{r}$ in $T_{G}$ gets sent to $c_{2r}$ by $\stand_{C}$ and does not affect the mapping of 
the other letters in the tableau. Hence, $S_{G} = \stand_{C}(T_{G})$.
\end{proof}

%%%%%%%%%%%%%%%%%%%%%%%%%%%%%%%%%%%%%%%%%%%%%%%%%%%%%%%%%%%%%%%%%
\section{Characterization of the shape of graphs}
\label{section.shape}

While answering Problem~\ref{problem.shape} in full generality seems far-achieving, we determine necessary and sufficient conditions of a 
hook graph in this section.

%%%%%%%%%%%%%%%%%%%%%%%%%%%%%%%%%%%%%%%%%%%%%%%%%%%%%%%%%%%%%%%%%
\subsection{Trees}

We begin by establishing a necessary condition for a connected graph to be of hook shape.

\begin{proposition}
\label{proposition.connected}
Let $G$ be a simple graph and let $k$ be the number of connected components of $G$ that contain at least one edge. If $G$ contains $k$ 
edges $e_{1}, \ldots, e_{k}$ such that $G-\{e_{1}, \ldots, e_{k}\}$ has the same number of connected components as $G$, then $G$ is not a hook-graph.
\end{proposition}

\begin{proof}
Let $C_{1}, \ldots, C_{k}$ denote the $k$ connected components of $G$ that contain at least one edge. Denote by $n_{i}$ the number of vertices in 
$C_{i}$ and let $n = \sum_{i=1}^{k} n_{i}$. If each $C_{i}$ was minimally connected, it would contain $n_{i}-1$ edges. Thus, $G$ contains at least 
$\sum_{i=1}^{k} (n_{i}-1) = n-k$ edges.  The condition that ``$G$ contains $k$ edges $e_{1}, \ldots, e_{k}$ such that $G-\{e_{1}, \ldots, e_{k}\}$ 
has the same number of connected components as $G$" implies that $G$ has at least $n-k+k = n$ edges. Observe that the length of $d_{G}$ is 
precisely $n$. However, the length of the partition $(e, 1^{e})$ is at least $n+1$, where $e\geqslant n$ is the number of edges of $G$. This implies that 
the partition $(e, 1^{e})$ is not weakly greater than $\tilde{d}_{G}$ in dominance order. Thus, there are no semistandard Young tableaux of shape $(e, 1^{e})$ 
with weight $d_{G}$, and $T_{G}$ is not hook-shaped.
\end{proof}

Recall that an undirected graph is called a \defn{tree} if it is connected and does not contain any cycle.
Setting $k=1$ into Proposition~\ref{proposition.connected}, we obtain the following necessary condition for connected (excluding singletons) hook-graphs.

\begin{corollary} 
The only connected hook-graphs are trees.
\end{corollary}

\begin{remark}
\label{remark.tree}
A tree need not be a hook-graph always. For example, the tree whose Burge array is 
$\begin{bmatrix}
    2  & 4 & 4 \\
    1  & 3 & 2
\end{bmatrix}$ has the shape $(2,2,2)$.
\end{remark}

%%%%%%%%%%%%%%%%%%%%%%%%%%%%%%%%%%%%%%%%%%%%%%%%%%%%%%%%%%%%%%%%%
\subsection{Peak and valley condition}

We now introduce peak and valley conditions which characterize when a graph has hook shape.

\begin{definition}[Peak]
\label{def peak}
A simple graph $G$ with $\AG$ as in Definition~\ref{burge array} is said to have a \defn{peak} if there exist $1\leqslant i < j <k\leqslant r$ such that 
\begin{enumerate}
\item $b_{i} \leqslant b_{k}$, 
\item $j$ is the minimum index with $b_{k} < b_{j}$,
\item $a_{i} \leqslant b_{j}$.
\end{enumerate}
\end{definition}

\begin{definition}[Valley]
\label{def valley}
A simple graph $G$ with $\AG$ as in Definition~\ref{burge array} is said to have a \defn{valley} if there exist $1\leqslant i < j <k\leqslant r$ such that the
following conditions hold:
\begin{enumerate}
\item $b_j \leqslant b_k< a_j$,
\item $b_j < b_i.$
\end{enumerate}
\end{definition}

\begin{example}
\mbox{}
\begin{enumerate}
\item
The graph with Burge array
$\begin{bmatrix}
   2  & 4 & 4 \\
    1  & 3 & 2
\end{bmatrix}$ 
of Remark~\ref{remark.tree} has a peak with $i=1,j=2,k=3$,  but no valley.
\item
The graph $G$ with 
$\AG= \begin{bmatrix}
    4& 5 & 6& 7\\
    1& 3 & 5& 2
\end{bmatrix}$ does not have a peak as $b_1<b_4<b_2$ but $a_1>b_2$. Note that $b_1<b_4<b_3$, but $j=3$ is not the minimal $j$ satisfying this condition.
Also, $\AG$ does not have a valley.
\item
The graph considered in Example \ref{eg nonhook} has both a peak and a valley.
\end{enumerate}
\end{example}

We refer to a Burge array as being \defn{PV-free} if the Burge array does not contain a peak or a valley.

\begin{theorem}
\label{thm.hook_char}
Let $G$ be a simple graph on $[n]$. 
The graph $G$ is a hook-graph if and only if its corresponding Burge array is PV-free.
\end{theorem}

\begin{proof} 
We prove the equivalent statement that the shape of $G$ is \emph{non-hook} if and only if $G$ has either a peak or a valley.

\smallskip
\noindent
\textbf{Proof of forward direction $\Rightarrow$:} 

\noindent
Let $T_G$ be the tableau of non-hook shape associated with the Burge array $\AG=$ 
$
\begin{bmatrix}
    a_1       & a_2 & \dots & a_{r} \\
    b_1       & b_2 & \dots & b_r
\end{bmatrix}$.
Let $T_\ell$ denote the tableau corresponding to the sub-array consisting of the first $\ell$ columns of $\AG$.
Choose $k$ minimal such that the shape of $T_k$ is non-hook, that is, let $k$ be the column that creates the cells $(2,2)$ and $(3,2)$ when applying
the Burge algorithm. We claim that there exist  $j_1,j_2$ with $1\leqslant j_1<j_2< k$ such that $
\begin{bmatrix}
    a_{j_1}       & a_{j_2} & a_k \\
    b_{j_1}       & b_{j_2} & b_k
\end{bmatrix}$
is either a peak or a valley.

Let $x$ be the first row entry of $T_{k-1}$ that is bumped by $b_k$ in $k$-th step of the Burge algorithm. Let $y$ be the entry in position $(2,1)$ 
(first entry of 2nd row) of $T_{k-1}$, see Figure \ref{T_k-1 and T_k}. Note that $x>b_k$ and $y\leqslant x$.
\begin{figure}[h]
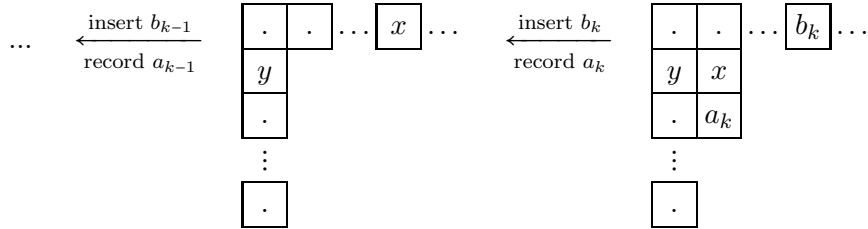

\label{T_k-1 and T_k}
...
\quad
$\xleftarrow[\text{record } a_{k-1}]{\text{insert } b_{k-1}}$
\quad
\begin{ytableau}
. & . & \none[\dots] & x & \none[\dots]\\
y \\
. \\
\none[\vdots]\\
.
\end{ytableau}
\quad 
$\xleftarrow[\text{record } a_k]{\text{insert } b_k}$
\quad
\begin{ytableau}
. & . & \none[\dots] & b_k & \none[\dots]\\
y & x\\
. & a_k\\
\none[\vdots]\\
.
\end{ytableau}
\caption{Tableau after $(k-1)$ insertions (left) and tableau after $k$ insertions (right)
\label{T_k-1 and T_k}}
\end{figure}
There are two different cases depending on whether $x$ is an inserted letter or a recorded letter.

\smallskip
\noindent
\textbf{Case 1:} Let $x$ be the inserted letter $b_j$ for some $1<j\leqslant k-1$. 

\noindent
\textbf{Subcase A:} If $y$ is a recorded letter, then the only possibility is $y=a_1$. This implies that position $(1,1)$ of $T_{k-1}$ is $b_1$ and the arm 
of $T_{k-1}$ consists of $b_{2} b_{3} \ldots b_{k-1}$. So we have that $j$ is the smallest index between $1$ and $k$ satisfying $b_1\leqslant b_k<b_j$. 
In addition, $a_1 =y \leqslant x= b_j$. Hence 
 $\begin{bmatrix}
    a_{1}       & a_{j} & a_k \\
    b_{1}       & b_{j} & b_k
\end{bmatrix}$
is a peak. 

\noindent
\textbf{Subcase B:} 
Assume $y$ is an inserted letter, say $b_i$ for some $i$. Let $b_m$ be the value which bumped $b_i$ from the first row. So $b_m<b_i \leqslant x=b_j$. 
Also we must have $b_m\leqslant b_k$ as otherwise $b_k< b_m < x$ and in this case $b_k$ would bump an entry strictly smaller than $x=b_j$. Now two 
cases arise: $m<j$ and $m>j$. When $m>j$, we have $a_j<a_m$. We use the relations $b_m\leqslant b_k<x=b_j<a_j<a_m$  and $b_m<b_i$ to show that 
$\begin{bmatrix}
    a_{i}       & a_{m} & a_k \\
    b_{i}       & b_{m} & b_k
\end{bmatrix}$ 
is a valley. When $m<j$, we claim that 
$\begin{bmatrix}
    a_{m}       & a_{j} & a_k \\
    b_{m}       & b_{j} & b_k
\end{bmatrix}$
is a peak. First observe that $a_m\leqslant b_j$. Since $m<j$, $x=b_j$ is inserted after $b_m$. If $a_m>b_j$ then $b_j$ would bump an entry weakly 
larger than $y$ and the shape of $T_{k-1}$ becomes non-hook. This contradicts the minimality of $k$. Hence $a_m\leqslant b_j$. Also note that $j$ is 
the smallest index between $m$ and $k$ satisfying $b_m\leqslant b_k<b_j$ as otherwise this would contradict the choice of $y$ or $x$. Hence the claim is true.

\smallskip
\noindent
\textbf{Case 2:} Let $x$ be a recorded letter say $a_j$ for $1<j\leqslant k-1$. This implies $y=b_i$ for some $i$ as $y$ can no longer be the recording 
letter $a_1$.

\noindent
\textbf{Subcase A: $i\leqslant j$.} Let $b_m$ be the value which bumps $b_i$ from the first row. Then we must have $m\geqslant j$. If $m<j$, then 
$b_i$ would be placed at position $(2,1)$ of $T_m$. Since $x=a_j$ is placed in the first row, $b_j$ must bump a letter from the first row. That letter 
would bump $b_i$ from position $(2,1)$ of $T_j$. This contradicts that $y=b_i$. Now we have the relations $b_m\leqslant b_k<x=a_j\leqslant a_m$ 
and $b_m<b_i$. Hence 
$\begin{bmatrix}
    a_{i}       & a_{m} & a_k \\
    b_{i}       & b_{m} & b_k
\end{bmatrix}$
is a valley.

\noindent
\textbf{Subcase B: $i>j$.} Since $x=a_j$ is placed in the first row, $b_j$ must bump a letter from the first row. This letter must be an inserted letter 
since $a_1$ is already placed in the first column and hence a bumped recorded letter would create a non-hook shape. This would contradict the minimality 
of $k$. Let $b_t$ with $t<j$ be the value bumped by
$b_j$. So $b_t>b_j$. Let $b_m$ denote the value that bumps $y=b_i$. Here $m>i$ but $b_m<b_i$. Observe that $b_i<b_t$, as $b_i$ and $b_t$ are in the 
first column. We also have $b_m\leqslant b_k$ otherwise $b_k$ would bump something smaller than $x=a_j$. Moreover, $b_k<a_j\leqslant a_i\leqslant a_m$. 
Therefore $b_m\leqslant b_k<a_m$ and $b_t>b_m$, which shows that
$\begin{bmatrix}
    a_{t}       & a_{m} & a_k \\
    b_{t}       & b_{m} & b_k
\end{bmatrix}$
is a valley.

\smallskip
\noindent
\textbf{Proof of backward direction $\Leftarrow$:} 

\noindent
Now we shall prove that the shape of $T_G$ is non-hook if $G$ has either a valley or a peak. 
We start by proving the result for the case of a valley.

\smallskip
Suppose $\AG$ contains a \textbf{valley} $
\begin{bmatrix}
    a_i      & a_j & a_k \\
    b_i      & b_j & b_k
\end{bmatrix}$. We may assume that there is no peak or valley in the first $k-1$ columns of $\AG$ and that $T_{k-1}$ has hook shape.
We consider two cases.

\smallskip

\noindent
\textbf{Case 1:} Assume that $b_i$ is present in the first row of $T_{j-1}$. Since $b_i>b_j$, $b_j$ bumps an element, say $x$, from the first row of $T_{j-1}$. 
So $b_j<x\leqslant b_i$. Since $b_j$ bumps a letter in the first row and since by assumption $T_j$ has hook shape, it follows that $a_j$ must be in the first row 
of $T_j$. Note that $a_j$ must also be in the first row of $T_{k-1}$. This is because if an element say $b_\ell$ (insertion letter) bumps $a_j$ from the 
first row before the insertion of $b_k$, then $T_\ell$ has non-hook shape since $a_j\geqslant a_1$, which is greater or equal to the letter in cell $(2,1)$.

Since we have $b_j\leqslant b_k<a_j$, $b_k$ bumps an element (which lies in the first row of $T_{k-1}$),
say $y$. Then $b_j<y\leqslant a_j$. We claim that $x\leqslant y$. If $y<x$, note that $y$ was not in $T_{j-1}$ when $b_j$ was inserted, since otherwise $b_j$ 
would have bumped $y$. So $b_j<y<x \leqslant b_i<a_i<a_j$. This shows that  $y$ must be an inserted letter, say $b_\ell$ with $j<\ell<k$, and we have a valley 
formed by columns $i,j,\ell$. This is a contradiction to our assumption. Therefore we must have $x\leqslant y$. In this case the shape of $T_k$ 
becomes non-hook, because the entry in cell $(2,1)$ of $T_{k-1}$ is less than or equal to $x$ and so $y$ must be placed in the $(2,2)$ position of $T_k$.

\smallskip

\noindent
\textbf{Case 2:} Assume that $b_i$ is not present in the first row of $T_{j-1}$. Let $x$ be the element in $T_{j-1}$ in the position where $b_i$ was originally inserted. Since $x<b_i$ 
and $x$ is inserted after step $i$, $x$ must be an inserted letter, say $b_\ell$ with $i<\ell<j$. We claim that $b_\ell>b_j$. If not, we have $b_\ell\leqslant b_j$
and 
$\begin{bmatrix}
    a_i      & a_\ell & a_j \\
    b_i      & b_\ell & b_j
\end{bmatrix}$ 
is a valley since $b_j <b_i<a_i\leqslant a_\ell$. This is a contradiction. Hence $b_\ell>b_j$. This implies that $b_j$ must bump an element from the first row 
of $T_{j-1}$.  If $z$ is that element, then $b_j<z\leqslant b_\ell$. Also, since $b_j$ bumps an element, $a_j$ must be in the first row in $T_j$. By the same
arguments as in Case 1, $a_j$ must be in the first row in $T_{k-1}$. Since by the valley condition $b_k<a_j$, $b_k$ bumps an element in $T_{k-1}$.
Call this element $w$. We claim that $z \leqslant w$. If $w<z$, note that $w$ was not in $T_{j-1}$ since then $b_j\leqslant b_k<w<z$ and hence $b_j$ would
have bumped $w$ instead of $z$ in $T_{j-1}$. If $w$ is an inserted letter $b_m$ with $j<m<k$, then 
$\begin{bmatrix}
    a_i      & a_j & a_m \\
    b_i      & b_j & b_m
\end{bmatrix}$ 
forms a valley as $b_m=w<z\leqslant b_\ell <a_\ell\leqslant a_j$. This is a contradiction. If $w$ is a recorded letter, then $w\geqslant a_j 
\geqslant a_\ell> b_\ell \geqslant z$, contradicting the assumption that $w<z$. This proves $z \leqslant w$.
Hence $w$ must be placed in position $(2,2)$ of $T_k$. This implies that the shape of $T_k$ is non-hook.

Considering the above cases, we conclude that $T_G$ has non-hook shape when $\AG$ has a valley.

\smallskip

Suppose $\AG$ contains a \textbf{peak} $\begin{bmatrix}
    a_{i}       & a_{j} & a_k \\
    b_{i}       & b_{j} & b_k
\end{bmatrix}$. As before, we may assume that there is no valley or peak in first $k-1$ columns of $\AG$ and that $T_{k-1}$ has hook shape.

\smallskip
\noindent
\textbf{Claim:} $b_j$ cannot be bumped from the first row before the $k$-th step.

\noindent
\textit{Proof:}
Assume that $b_j$ is bumped from the first row before the $k$-th column is inserted. This would create a non-hook shape as 
$b_j \geqslant a_{i} \geqslant a_{1}$. This contradicts the fact that $T_{k-1}$ has hook shape. Hence we have proved the claim.

Since $b_j$ is in the first row of $T_{k-1}$ and $b_k<b_j$, $b_k$ must bump an element, say $z$. So $b_k<z\leqslant b_j$. 

\smallskip

\noindent
\textbf{Case 1:} Suppose $b_i$ is not present in the first row of $T_{k-1}$. Let $y$ be the entry in position $(2,1)$ of $T_{k-1}$. Then $y\leqslant b_i$. 
Also we have $b_i\leqslant b_k<z$. Since $y\leqslant z$, $z$ would be placed at position $(2,2)$ of $T_k$, hence the shape of $T_{k}$ becomes non-hook.

\smallskip

\noindent
\textbf{Case 2:} Assume $b_i$ is present in the first row of $T_{k-1}$.

\noindent
\textbf{Subcase A:} $b_i$ is  placed at the end of the first row of $T_{i-1}$. 
If $z=a_\ell$, then $z\geqslant a_{1}$ and $T_{k}$ has non-hook shape since the element in position $(2,1)$ in $T_{k-1}$ is smaller or equal to $a_1$.

Now consider the case when $z=b_\ell$ for some $\ell$. 
If $\ell\leqslant j$, consider the array 
$\begin{bmatrix}
    a_{i}       & a_{\ell}  & a_j   & a_k \\
    b_{i}       & b_{\ell}  & b_j   & b_k
\end{bmatrix}$.
Note that indeed $i< \ell$. Namely, $z = b_{\ell}$ implies $b_{i} < b_{\ell}$ since $b_{i} \leqslant b_{k}$ and $b_{k} < z$. This implies that $i < \ell$; otherwise, 
$b_{i}$ would bump a value from the first row of $T_{i-1}$. We have $b_i \leqslant b_k<z=b_\ell$. By definition of the peak 
$\begin{bmatrix}
    a_{i}       & a_{j} & a_k \\
    b_{i}       & b_{j} & b_k
\end{bmatrix}$, 
we have $a_i\leqslant b_\ell=z$. Note that by the definition of a peak $j=\ell$. Again, $y\leqslant z$ implies that shape of $T_k$ is non-hook.

Now let us consider the case when $\ell>j$. If $b_\ell=b_j$, then $b_{j}$ would be bumped by $b_{k}$ instead of $b_{\ell}$, which is a contradiction.
If $b_\ell<b_j$, then
$\begin{bmatrix}
    a_{i}       & a_{j} & a_\ell \\
    b_{i}       & b_{j} & b_\ell
\end{bmatrix}$ becomes a peak. Observe that $j$ is minimal due to the minimality condition of the peak $\begin{bmatrix}
    a_{i}       & a_{j} & a_k \\
    b_{i}       & b_{j} & b_k
\end{bmatrix}$ and the inequality $b_{k} < b_{\ell}$.
This contradicts our assumption that there should not be any other peak in the first $k$ columns of $\AG$.

\noindent
\textbf{Subcase B:} $b_i$ bumps an element from the first row of $T_{i-1}$ when inserted. Then $a_i$ is  placed in the first row of $T_i$. 
Let $u$ be the value bumped by $b_i$ when inserted to $T_{i-1}$ and let $z$ be the value bumped by $b_k$. So, $b_i<u$ and $b_k<z$. Notice that the 
entry in position $(2,1)$ of $T_{k-1}$ is less than or equal to $u$. So, if $u\leqslant z$, the shape of $T_k$ becomes non-hook.

On the contrary, suppose $u>z$. Then $z$ is not in the first row of $T_{i-1}$ since otherwise $b_i$ would bump $z$ instead of $u$ as $b_i\leqslant b_k<z<u$. 
This contradicts the definition of $u$. Hence $z$ is inserted after the $i$-th step. As $b_{i}$ bumps $u$ from the first row of $T_{i-1}$, $u$ is weakly less than 
$a_{i}$. As $z$ is inserted after the $i$-th step, it must be an inserted letter i.e. $z = b_{\ell}$ for some $i < \ell < k$. Now either $z<a_i$ or $z\geqslant a_i$. 

Assume that $z < a_{i}$. We prove that there does not exists an inserted letter $b_{t}$ such that $i < t < k$ and $b_{i} \leqslant b_{t} < a_{i}$. Assume that 
there exists such an inserted letter $b_{t}$, where $t$ is the smallest possible index. Observe that when $b_{t}$ is inserted it bumps an element strictly 
right of $b_{i}$ and weakly left of $a_{i}$ in the arm of $T_{t-1}$. Denote this element by $w$. From the minimality of $t$, we have $w$ is present in the first 
row of $T_{i-1}$ and is strictly to the right of $u$. Thus, $w$ is weakly larger than $u$, and when $w$ is bumped, it will be weakly larger than the value in 
the $(2,1)$ cell of $T_{t-1}$. This would contradict $T_{k-1}$ being hook shaped. Therefore no such inserted letter $b_{t}$ exists. Observe that $z$ satisfies 
the condition of $b_{t}$ namely $z = b_{\ell}$ where $i < \ell < k$ and $b_{i} \leqslant  b_{\ell} < a_{i}$. By our claim this is impossible. 

Next consider the case $z\geqslant a_i$. Note that $a_1 \geqslant u$ as both of them are placed in first column. Together with $u>z \geqslant a_i$ we 
get $a_1>a_i$, which is not possible. Therefore the case $u>z$ does not arise. 

Considering all the cases, we proved that the shape of $T_k$ is non-hook when $\AG$ has a peak.
This completes the proof.
\end{proof}

\begin{remark}
Theorem~\ref{thm.hook_char} is the analogue to the statement for RSK that the shape of a word $w$ under RSK is a single row if and only if $w$ is 
weakly increasing.
\end{remark}

\begin{remark}
In analogy with Schensted's result for the RSK insertion that the length of the longest increasing subsequence of a word $w$ gives the length
of the longest row in the Young tableaux under RSK, one might suspect that the longest PV-free subsequence of a Burge array
gives the size of the largest hook in $\sh(T_G)$. However, this is not true as the following counterexample shows. Take the graph $G$ with Burge array
\[
	\AG= \begin{bmatrix}
    	4& 8 & 8& 9& 9\\
    	1& 3 & 2& 5& 2
	\end{bmatrix}.
\]
The tableau under the Burge correspondence is
\[
	T_G=
	\begin{ytableau}
	1 & 2 & 2\\
	3 & 5 & 9\\
	4 & 8\\
	8 & 9
	\end{ytableau}\;,
\]
so that $\sh(T_G) = (3,3,2,2)$. However, the subsequence
\[
	\begin{bmatrix}
    	4& 8 & 9& 9\\
    	1& 3 & 5& 2
	\end{bmatrix}
\]
is PV-free and has hook shape $(4,1,1,1,1)$, which is larger than the biggest hook $(3,1,1,1)$ in $\sh(T_G)$.
\end{remark}

A \defn{star graph} is a graph where one vertex $i$ is connected to all other vertices by an edge and no other edges exist in the graph.

\begin{corollary} 
All star graphs are hook-graphs.
\end{corollary}

\begin{proof}
Assuming that the vertices of the star graph are labelled $1,2,\ldots,n$ and the vertex connected to all other vertices is vertex $i$, the Burge array is
\[
	\begin{bmatrix}
	i& i& \dots& i& i+1& i+2& \dots& n\\
	i-1& i-2& \dots& 1& i& i& \dots &i
	\end{bmatrix},
\]
which is PV-free.
\end{proof}

%%%%%%%%%%%%%%%%%%%%%%%%%%%%%%%%%%%%%%%%%%%%%%%%%%%%%%%%%%%%%%%%%
\section{Crystal Structure on hook-graphs}
\label{section.crystal}

We review the crystal structure on semistandard Young tableaux in Section~\ref{section.crystal review} and then define
the new crystal structure on hook-graphs in Section~\ref{section.crystal hook}.

%%%%%%%%%%%%%%%%%%%%%%%%%%%%%%%%%%%%%%%%%%%%%%%%%%%%%%%%%%%%%%%%%
\subsection{Review of crystal structure on semistandard Young tableaux}
\label{section.crystal review}

Crystal bases provide a combinatorial skeleton for $U_q(\mathfrak{g})$-representations, where $U_q(\mathfrak{g})$ is the quantum
group associated to the Lie algebra $\mathfrak{g}$. A \defn{Kashiwara crystal} is a nonempty set $B$ together with maps
\[
\begin{split}
	e_i, f_i &\colon B \to B \cup \{\emptyset\} \qquad \text{for $i\in I$},\\
	\mathsf{wt} &\colon B \to \Lambda,
\end{split}
\]
where $\Lambda$ is the weight lattice associated to the Lie algebra $\mathfrak{g}$ and $I$ is the index set of the Dynkin diagram for $\mathfrak{g}$.
These maps have to satisfy certain conditions (see~\cite[Definition 2.13]{BumpSchilling.2017}).

For simply-laced Lie algebras $\mathfrak{g}$, a \defn{Stembridge crystal} is a Kashiwara crystal for which the raising and lowering operators $e_i$ and 
$f_i$ satisfy certain local rules (see~\cite[Section 4.2]{BumpSchilling.2017}). Stembridge crystals are crystals corresponding to 
$U_q(\mathfrak{g})$-representations.

In this paper, we only consider crystals of type $A_{m-1}$. For type $A_{m-1}$ crystals we have $I=\{1,2,\ldots,m-1\}$. A particular model
for type $A_{m-1}$ crystals is given in terms of semistandard Young tableaux, that is, the set $B$ is the set of all semistandard Young tableaux on the
alphabet $[m]$. For further details, see~\cite[Section 3]{BumpSchilling.2017}. Here we review the crystal structure on semistandard Young tableaux.

\begin{definition} 
Let $T$ be a semistandard Young tableau. The \defn{reading word} of $T$ denoted by $R(T)$ is obtained by reading the entries within a row from left to 
right starting with the bottommost row. The \defn{$i$-th reading word} of $T$ denoted by $R_{i}(T)$ is the induced subword of $R(T)$ containing only the entries 
$i$ and $i+1$.
\end{definition}

\begin{definition}
Let $\omega$ be a word with length $n$ over the alphabet $[m]$. A \defn{Knuth move} on $\omega$ is one of the following transformations:
\begin{enumerate}
\item $\omega_{1} \ldots bca \ldots \omega_{n} \longrightarrow \omega_{1} \ldots bac \ldots \omega_{n}$ if $a < b \leqslant c$,
\item $\omega_{1} \ldots bac \ldots \omega_{n} \longrightarrow \omega_{1} \ldots bca \ldots \omega_{n}$ if $a < b \leqslant c$,
\item $\omega_{1} \ldots acb \ldots \omega_{n} \longrightarrow \omega_{1} \ldots cab \ldots \omega_{n}$ if $a \leqslant b < c$,
\item $\omega_{1} \ldots cab \ldots \omega_{n} \longrightarrow \omega_{1} \ldots acb \ldots \omega_{n}$ if $a \leqslant b < c$.
\end{enumerate}
Two words $\omega$ and $\nu$ are said to be \defn{Knuth equivalent} if they differ by a sequence of Knuth moves.
\end{definition}

It is well-known that $\omega$ and $\nu$ are Knuth equivalent if and only if $P(\omega) = P(\nu)$, that is, their insertion tableaux under Schensted 
insertion are equal. In addition, it is also known that the crystal operators $f_{i}$ and $e_{i}$ on words preserve Knuth equivalence. Thus, in order to 
define the crystal operators on semistandard Young tableaux it suffices to look at their reading words.

\begin{definition} 
\label{definition.crystal}
Let $T$ be a semistandard Young tableau in $\Tab_m(\lambda)$. Assign a `)' to every $i$ in $R_{i}(T)$ and a `(' to every $i+1$ in $R_{i}(T)$. 
Successively pair every `(' that is directly left of a `)' which we call an \defn{$i$-pair} and remove the $i$-paired terms. Continue this process until no more 
terms can be $i$-paired. 

The \defn{lowering operator $f_{i}$} for $1\leqslant i<m$ acts on $T$ as follows:
\begin{enumerate}
\item If there are no unpaired `)' terms left, then $f_{i}$ annihilates $T$.
\item Otherwise locate the $i$ in $T$ corresponding to the rightmost unpaired `)' of $R_{i}(T)$ and replace it with an $i+1$.
\end{enumerate}

The \defn{raising operator $e_{i}$} for $1\leqslant i<m$ acts on $T$ as follows:
\begin{enumerate}
\item If there are no unpaired `(' terms left, then $e_{i}$ annihilates $T$.
\item Otherwise locate the $i+1$ in $T$ corresponding to the leftmost unpaired `(' of $R_{i}(T)$ and replace it with an $i$.
\end{enumerate}
The \defn{weight} $\mathsf{wt}(T)=(a_1,a_2,\ldots,a_m)$ is an $m$-tuple such that $a_i$ is the number of letters $i$ in $T$.
\end{definition}

The crystal lowering  and raising operators $f_i$ and $e_i$ for $1\leqslant i <m$ together with the weight function $\mathsf{wt}$
define a crystal structure on $\Tab_m(\lambda)$. The vertices of the 
crystal are the elements in $\Tab_m(\lambda)$ for a fixed partition $\lambda$. There is an edge labelled $i$ from $T\in \Tab_m(\lambda)$ to 
$T'\in \Tab_m(\lambda)$ if $f_i(T)=T'$. Note that $f_i$ and $e_i$ are partial inverses, that is, if $f_i(T) = T'$ then $e_i(T')=T$ and vice versa.

%%%%%%%%%%%%%%%%%%%%%%%%%%%%%%%%%%%%%%%%%%%%%%%%%%%%%%%%%%%%%%%%%
\subsection{Crystal structure on hook-graphs}
\label{section.crystal hook}

In this section, we assume that $G$ is a hook-graph or equivalently by Theorem~\ref{thm.hook_char} that $\mathcal{A}_G$ is a $PV$-free Burge array.
 
\begin{definition} 
\label{definition.reading word AG}
The $i$-th reading word of $\mathcal{A}_G$, denoted by $\tilde{R}_{i}(\mathcal{A}_G)$, is obtained by the following algorithm: 
\begin{enumerate}
\item Let $a_{k}$ denote the leftmost $i+1$ in the top row of $\mathcal{A}_G$. If $k = 1$ or $b_{k-1} \leqslant b_{k}$, then let $a_{k}$ be the first 
letter of $\tilde{R}_{i}(\mathcal{A}_G)$.
\item Read all other $i$'s and $(i+1)$'s in $\mathcal{A}_G$ from left to right while appending the corresponding value to $\tilde{R}_{i}(\mathcal{A}_G)$.
\end{enumerate}
\end{definition}

\begin{example}
Let
$\AG=\begin{bmatrix}
   3&3&4\\
   2&1&3
\end{bmatrix}$. Then $\tilde{R}_1=21$, $\tilde{R}_2=3233$, and $\tilde{R}_3=4333$. For
$\AG=\begin{bmatrix}
   3&4\\
   2&1
\end{bmatrix}$, we have $\tilde{R}_3=34$.
\end{example}

\begin{remark}
Note that except for the column in $\AG$ containing $a_k$ as in (1) of Definition~\ref{definition.reading word AG}, each column of $\AG$ contains
either $i$ or $i+1$, but not both. Hence the algorithm to construct the reading word in Definition~\ref{definition.reading word AG} is well-defined.
Indeed, if $\AG$ contains the column $\begin{bmatrix} i+1\\ i \end{bmatrix}$, then it must be the leftmost column containing $i+1$ in the top row
since by the definition of a Burge array the bottom row is decreasing for equal top row elements. Hence $a_k$ is this leftmost $i+1$ and either $k=1$ 
or $a_{k-1}\leqslant i$ and $b_{k-1}<i=b_k$, so that $a_k$ is chosen as the first letter of $\tilde{R}_i(\AG)$.
\end{remark}

\begin{definition} 
\label{definition.crystal operators}
Assign a `)' to every $i$ in $\tilde{R}_{i}(\mathcal{A}_G)$ and a `(' to every $i+1$ in $\tilde{R}_{i}(\mathcal{A}_G)$. Successively pair every `(' that 
is directly left of a `)', called an \defn{$i$-pair}, and remove the paired terms. Continue this process until no more terms can be paired.

The operator $\tilde{f}_{i}$ acts on $\mathcal{A}_{G}$ as follows:
\begin{enumerate}
\item If there are no unpaired `)' terms left, then $\tilde{f}_{i}$ annihilates $\mathcal{A}_{G}$ denoted by $\tilde{f}_{i}(\mathcal{A}_G) = 0$.
\item Otherwise locate the $i$ in $\mathcal{A}_{G}$ corresponding to the rightmost unpaired `)' and denote it by $x$.
	\begin{enumerate}
	\item If there is no $i+1$ in the same column as $x$, then $\tilde{f}_{i}$ changes $x$ in $\mathcal{A}_{G}$ to an $i+1$.
	\item If there is an $i+1$ in the same column as $x$, then $x$ is on the bottom row of  $\mathcal{A}_{G}$. Let $k$ be the index such that 
	$b_{k} = x$. Let $\ell$ be the smallest index such that $b_{\ell} \leqslant b_{\ell +1} \leqslant \cdots \leqslant b_{k-1}$. Let $\ell \leqslant m \leqslant k$ 
	be the largest index such that $b_{m} < a_{\ell}$. In the top row of $\mathcal{A}_{G}$ replace $a_{k-1}$ with an $i+1$ and replace $a_{s}$ with 
	$a_{s+1}$ for $\ell \leqslant s \leqslant k-2$. In the bottom row of $\mathcal{A}_{G}$ replace $b_{m}$ with $a_{\ell}$ and replace $b_{k}$ with $b_{m}$.
	(Remark: Observe that in this case $k \not = 1$ as otherwise the $i+1$ and $i$ in this column would form an $i$-pair in $\tilde{R}_{i}$. Thus, this 
	operation is well-defined.) This procedure is illustrated in Figure~\ref{figure.fi action}.
	\end{enumerate}
\end{enumerate}

\begin{figure}[h]
\begin{displaymath}
\begin{bmatrix}
    \dots \,\, a_{\ell}  & a_{\ell +1} & \dots & a_{m-1} & a_{m} & a_{m+1} & \dots & a_{k-1} & i+1 &\dots\\
    \dots \,\, b_{\ell}  & b_{\ell +1} & \dots & b_{m-1} & b_{m} & b_{m+1} & \dots & b_{k-1} & i   &\dots
\end{bmatrix}
\end{displaymath}
$\downarrow{\tilde{f}_i}$
\begin{displaymath}
\begin{bmatrix}
    \dots \,\, a_{\ell +1}  & a_{\ell +2} & \dots & a_{m}    & a_{m+1}  & a_{m+2} & \dots & i+1     & i+1 &\dots \\
    \dots \,\, b_{\ell}     & b_{\ell +1} & \dots & b_{m-1}  & a_{\ell} & b_{m+1} & \dots & b_{k-1} & b_{m} & \dots
\end{bmatrix}
\end{displaymath}
\caption{Action of $\tilde{f}_i$ in Case (2b).
\label{figure.fi action}}
\end{figure}

The operator $\tilde{e}_{i}$ acts on $\mathcal{A}_{G}$ as follows:
\begin{enumerate}
\item If there are no unpaired `(' terms left, then $\tilde{e}_{i}$ annihilates $\mathcal{A}_{G}$ denoted by $\tilde{e}_{i}(\mathcal{A}_G) = 0$.
\item Otherwise locate the $i+1$ in $\mathcal{A}_{G}$ corresponding to the leftmost unpaired `(' and denote it by $x$.
	\begin{enumerate}
	\item If $x$ is in the top row of $\mathcal{A}_{G}$ and there is an $i+1$ directly to the left of it, then let $k$ be the index such that $a_{k} = x$. 
	Let $\ell$ be the smallest index such that $b_{\ell} \leqslant  b_{\ell +1} \leqslant \cdots \leqslant b_{k-1}$. Let $\ell \leqslant m \leqslant k$ be the smallest 
	index such that $b_{k} < b_{m}$. In the top row of $\mathcal{A}_{G}$ replace $a_{\ell}$ with $b_{m}$ and replace $a_{s}$ with $a_{s-1}$ for 
	$\ell+1 \leqslant s \leqslant k-1$. In the bottom row of $\mathcal{A}_{G}$ replace $b_{m}$ with $b_{k}$ and replace $b_{k}$ with an $i$.
	\item Otherwise, $\tilde{e}_{i}$ changes $x$ in $\mathcal{A}_{G}$ to an $i$.
	\end{enumerate}
\end{enumerate}
\end{definition}

\begin{figure}[t]
\begin{subfigure}[b]{.4\linewidth}
\scalebox{0.7}{
\begin{tikzpicture}[>=latex,line join=bevel,]
\node (node_0) at (93.0bp,280.0bp) [draw,draw=none] {$\left(\begin{array}{rr}2 & 4 \\1 & 3\end{array}\right)$};
  \node (node_1) at (62.0bp,104.0bp) [draw,draw=none] {$\left(\begin{array}{rr}3 & 4 \\2 & 3\end{array}\right)$};
  \node (node_2) at (93.0bp,368.0bp) [draw,draw=none] {$\left(\begin{array}{rr}2 & 4 \\1 & 2\end{array}\right)$};
  \node (node_3) at (25.0bp,192.0bp) [draw,draw=none] {$\left(\begin{array}{rr}3 & 4 \\2 & 2\end{array}\right)$};
  \node (node_4) at (93.0bp,192.0bp) [draw,draw=none] {$\left(\begin{array}{rr}3 & 4 \\1 & 3\end{array}\right)$};
  \node (node_5) at (62.0bp,456.0bp) [draw,draw=none] {$\left(\begin{array}{rr}2 & 4 \\1 & 1\end{array}\right)$};
  \node (node_6) at (161.0bp,280.0bp) [draw,draw=none] {$\left(\begin{array}{rr}3 & 4 \\2 & 1\end{array}\right)$};
  \node (node_7) at (130.0bp,456.0bp) [draw,draw=none] {$\left(\begin{array}{rr}2 & 3 \\1 & 2\end{array}\right)$};
  \node (node_8) at (25.0bp,280.0bp) [draw,draw=none] {$\left(\begin{array}{rr}3 & 4 \\1 & 2\end{array}\right)$};
  \node (node_9) at (25.0bp,368.0bp) [draw,draw=none] {$\left(\begin{array}{rr}3 & 4 \\1 & 1\end{array}\right)$};
  \node (node_10) at (96.0bp,544.0bp) [draw,draw=none] {$\left(\begin{array}{rr}2 & 3 \\1 & 1\end{array}\right)$};
  \node (node_11) at (161.0bp,368.0bp) [draw,draw=none] {$\left(\begin{array}{rr}3 & 3 \\2 & 1\end{array}\right)$};
  \node (node_12) at (161.0bp,192.0bp) [draw,draw=none] {$\left(\begin{array}{rr}4 & 4 \\2 & 1\end{array}\right)$};
  \node (node_13) at (130.0bp,104.0bp) [draw,draw=none] {$\left(\begin{array}{rr}4 & 4 \\3 & 1\end{array}\right)$};
  \node (node_14) at (99.0bp,16.0bp) [draw,draw=none] {$\left(\begin{array}{rr}4 & 4 \\3 & 2\end{array}\right)$};
  \draw [red,->] (node_0) ..controls (93.0bp,251.38bp) and (93.0bp,232.88bp)  .. (node_4);
  \definecolor{strokecol}{rgb}{0.0,0.0,0.0};
  \pgfsetstrokecolor{strokecol}
  \draw (102.0bp,236.0bp) node {$2$};
  \draw [green,->] (node_1) ..controls (74.139bp,75.128bp) and (82.116bp,56.156bp)  .. (node_14);
  \draw (92.0bp,60.0bp) node {$3$};
  \draw [red,->] (node_2) ..controls (93.0bp,339.38bp) and (93.0bp,320.88bp)  .. (node_0);
  \draw (102.0bp,324.0bp) node {$2$};
  \draw [red,->] (node_3) ..controls (37.139bp,163.13bp) and (45.116bp,144.16bp)  .. (node_1);
  \draw (56.0bp,148.0bp) node {$2$};
  \draw [blue,->] (node_4) ..controls (82.998bp,170.57bp) and (80.223bp,164.06bp)  .. (78.0bp,158.0bp) .. controls (74.622bp,148.8bp) and (71.436bp,138.52bp)  .. (node_1);
  \draw (87.0bp,148.0bp) node {$1$};
  \draw [green,->] (node_4) ..controls (105.14bp,163.13bp) and (113.12bp,144.16bp)  .. (node_13);
  \draw (123.0bp,148.0bp) node {$3$};
  \draw [blue,->] (node_5) ..controls (68.334bp,429.3bp) and (72.328bp,414.59bp)  .. (77.0bp,402.0bp) .. controls (78.1bp,399.04bp) and (79.356bp,395.98bp)  .. (node_2);
  \draw (86.0bp,412.0bp) node {$1$};
  \draw [red,->] (node_5) ..controls (49.861bp,427.13bp) and (41.884bp,408.16bp)  .. (node_9);
  \draw (56.0bp,412.0bp) node {$2$};
  \draw [green,->] (node_6) ..controls (161.0bp,251.38bp) and (161.0bp,232.88bp)  .. (node_12);
  \draw (170.0bp,236.0bp) node {$3$};
  \draw [green,->] (node_7) ..controls (117.86bp,427.13bp) and (109.88bp,408.16bp)  .. (node_2);
  \draw (123.0bp,412.0bp) node {$3$};
  \draw [red,->] (node_7) ..controls (140.13bp,427.25bp) and (146.73bp,408.52bp)  .. (node_11);
  \draw (157.0bp,412.0bp) node {$2$};
  \draw [blue,->] (node_8) ..controls (25.0bp,251.38bp) and (25.0bp,232.88bp)  .. (node_3);
  \draw (34.0bp,236.0bp) node {$1$};
  \draw [blue,->] (node_9) ..controls (25.0bp,339.38bp) and (25.0bp,320.88bp)  .. (node_8);
  \draw (34.0bp,324.0bp) node {$1$};
  \draw [green,->] (node_10) ..controls (84.845bp,515.13bp) and (77.515bp,496.16bp)  .. (node_5);
  \draw (91.0bp,500.0bp) node {$3$};
  \draw [blue,->] (node_10) ..controls (107.16bp,515.13bp) and (114.49bp,496.16bp)  .. (node_7);
  \draw (125.0bp,500.0bp) node {$1$};
  \draw [green,->] (node_11) ..controls (161.0bp,339.38bp) and (161.0bp,320.88bp)  .. (node_6);
  \draw (170.0bp,324.0bp) node {$3$};
  \draw [red,->] (node_12) ..controls (150.87bp,163.25bp) and (144.27bp,144.52bp)  .. (node_13);
  \draw (157.0bp,148.0bp) node {$2$};
  \draw [blue,->] (node_13) ..controls (119.87bp,75.253bp) and (113.27bp,56.517bp)  .. (node_14);
  \draw (126.0bp,60.0bp) node {$1$};
\end{tikzpicture}
}
\caption{Crystal of Burge arrays of shape $(2,1,1)$ with letters in $\{1,2,3,4\}$.}
\end{subfigure}
\begin{subfigure}[b]{.4\linewidth}
\scalebox{0.7}{
\begin{tikzpicture}[>=latex,line join=bevel,]
\node (node_0) at (75.0bp,544.0bp) [draw,draw=none] {$\left(\begin{array}{rrr}2 & 3 & 4 \\1 & 1 & 1\end{array}\right)$};
  \node (node_1) at (33.0bp,280.0bp) [draw,draw=none] {$\left(\begin{array}{rrr}3 & 4 & 4 \\1 & 2 & 1\end{array}\right)$};
  \node (node_2) at (33.0bp,368.0bp) [draw,draw=none] {$\left(\begin{array}{rrr}2 & 3 & 4 \\1 & 1 & 3\end{array}\right)$};
  \node (node_3) at (117.0bp,192.0bp) [draw,draw=none] {$\left(\begin{array}{rrr}3 & 3 & 4 \\2 & 1 & 3\end{array}\right)$};
  \node (node_4) at (33.0bp,192.0bp) [draw,draw=none] {$\left(\begin{array}{rrr}3 & 4 & 4 \\2 & 2 & 1\end{array}\right)$};
  \node (node_5) at (117.0bp,280.0bp) [draw,draw=none] {$\left(\begin{array}{rrr}2 & 3 & 4 \\1 & 2 & 3\end{array}\right)$};
  \node (node_6) at (80.0bp,16.0bp) [draw,draw=none] {$\left(\begin{array}{rrr}4 & 4 & 4 \\3 & 2 & 1\end{array}\right)$};
  \node (node_7) at (75.0bp,456.0bp) [draw,draw=none] {$\left(\begin{array}{rrr}2 & 3 & 4 \\1 & 1 & 2\end{array}\right)$};
  \node (node_8) at (80.0bp,104.0bp) [draw,draw=none] {$\left(\begin{array}{rrr}3 & 4 & 4 \\2 & 3 & 1\end{array}\right)$};
  \node (node_9) at (117.0bp,368.0bp) [draw,draw=none] {$\left(\begin{array}{rrr}2 & 3 & 4 \\1 & 2 & 2\end{array}\right)$};
  \draw [blue,->] (node_0) ..controls (75.0bp,515.38bp) and (75.0bp,496.88bp)  .. (node_7);
  \definecolor{strokecol}{rgb}{0.0,0.0,0.0};
  \pgfsetstrokecolor{strokecol}
  \draw (84.0bp,500.0bp) node {$1$};
  \draw [blue,->] (node_1) ..controls (33.0bp,251.38bp) and (33.0bp,232.88bp)  .. (node_4);
  \draw (42.0bp,236.0bp) node {$1$};
  \draw [green,->] (node_2) ..controls (33.0bp,339.38bp) and (33.0bp,320.88bp)  .. (node_1);
  \draw (42.0bp,324.0bp) node {$3$};
  \draw [blue,->] (node_2) ..controls (61.038bp,338.63bp) and (80.057bp,318.7bp)  .. (node_5);
  \draw (91.0bp,324.0bp) node {$1$};
  \draw [green,->] (node_3) ..controls (104.86bp,163.13bp) and (96.884bp,144.16bp)  .. (node_8);
  \draw (110.0bp,148.0bp) node {$3$};
  \draw [red,->] (node_4) ..controls (48.487bp,163.0bp) and (58.746bp,143.79bp)  .. (node_8);
  \draw (70.0bp,148.0bp) node {$2$};
  \draw [red,->] (node_5) ..controls (117.0bp,251.38bp) and (117.0bp,232.88bp)  .. (node_3);
  \draw (126.0bp,236.0bp) node {$2$};
  \draw [green,->] (node_5) ..controls (88.962bp,250.63bp) and (69.943bp,230.7bp)  .. (node_4);
  \draw (91.0bp,236.0bp) node {$3$};
  \draw [red,->] (node_7) ..controls (61.22bp,427.13bp) and (52.165bp,408.16bp)  .. (node_2);
  \draw (67.0bp,412.0bp) node {$2$};
  \draw [blue,->] (node_7) ..controls (88.78bp,427.13bp) and (97.835bp,408.16bp)  .. (node_9);
  \draw (108.0bp,412.0bp) node {$1$};
  \draw [green,->] (node_8) ..controls (80.0bp,75.378bp) and (80.0bp,56.877bp)  .. (node_6);
  \draw (89.0bp,60.0bp) node {$3$};
  \draw [red,->] (node_9) ..controls (117.0bp,339.38bp) and (117.0bp,320.88bp)  .. (node_5);
  \draw (126.0bp,324.0bp) node {$2$};
\end{tikzpicture}
}
\caption{Crystal of Burge arrays of shape $(3,1,1,1)$ with letters in $\{1,2,3,4\}$.}
\end{subfigure}
\caption{Examples of crystals on Burge arrays. \label{figure.crystals}}
\end{figure}
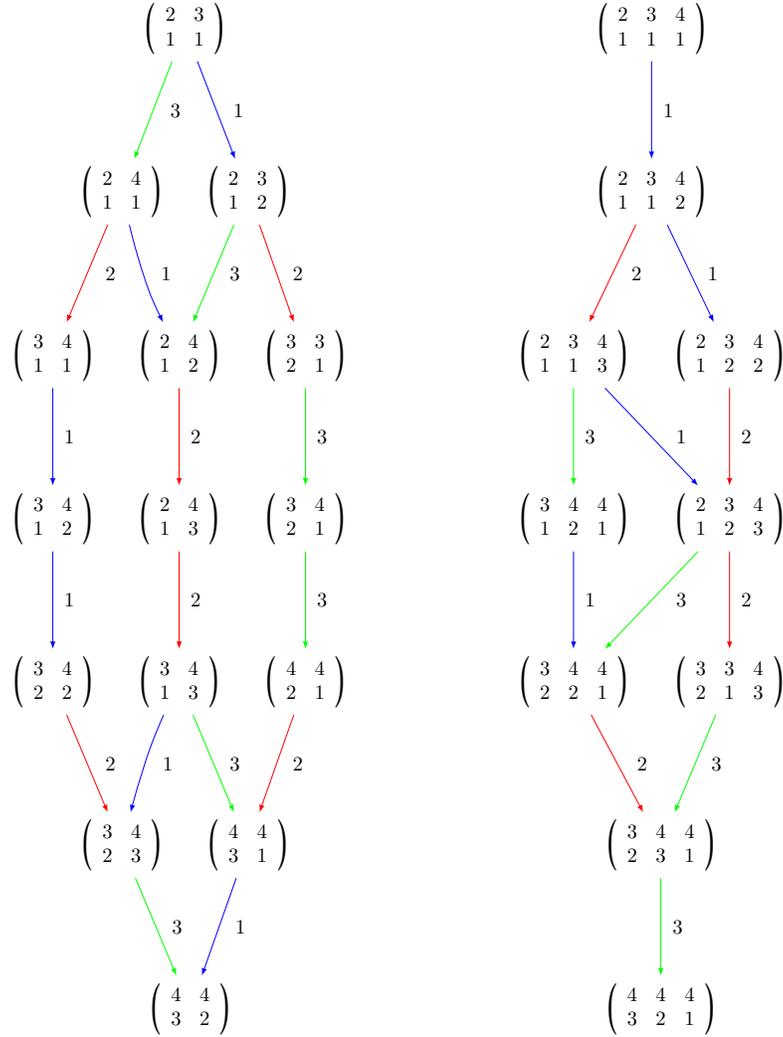

\begin{example}
Examples of crystals on Burge arrays are given in Figure~\ref{figure.crystals}. To illustrate the crystal operators $\tilde{f}_i$ of 
Definition~\ref{definition.crystal operators}, consider $\tilde{f}_2$ on 
$\AG= \begin{bmatrix}
   2&3&4\\
   1&2&3
\end{bmatrix}$. In this case $\tilde{R}_2(\AG) = 3223$, $x$ is the $2$ in column two of $\AG$, $k=2$, $\ell=1$, and $m=1$. We obtain
$\tilde{f}_2(\AG)= \begin{bmatrix}
   3&3&4\\
   2&1&3
\end{bmatrix}$. For $\tilde{f}_3$ on $\AG$ we have $x=3$, $k=3$, $\ell=1$, $m=1$, and $\tilde{f}_3(\AG) = \begin{bmatrix}
   3&4&4\\
   2&2&1
\end{bmatrix}$.
\end{example}

The fact that a non-annihilated $\tilde{f}_{i}(\mathcal{A}_G)$ (resp. $\tilde{e}_{i}(\mathcal{A}_G)$) is a $PV$-free Burge array or even a valid 
Burge array will be shown as a consequence of Proposition ~\ref{prop.hook_intertwine}. 

\begin{lemma} \label{lemma.i_pair}
$\tilde{R}_{i}(\mathcal{A}_G)$ has at most one $i$-pair.
\end{lemma}

\begin{proof}
Assume that $\tilde{R}_{i}(\mathcal{A}_G)$ has at least two $i$-pairs. This implies that $\mathcal{A}_G$ contains at least two $(i+1)$'s. Let $j$ be the 
index of the column containing the leftmost $i+1$ and $k$ be the index of the column containing the second leftmost $i+1$. We break into cases 
based on the position of the $(i+1)$'s in columns $j$ and $k$.

\smallskip

\noindent
\textbf{Case 1:} Assume that $a_{j} = i+1$ and $a_{k} = i+1$.
This implies that $a_{k}$ is the second leftmost $i+1$ in $\tilde{R}_{i}(\mathcal{A}_G)$. Since $\mathcal{A}_G$ is assumed to have at least two $i$-pairs, 
there must exist $\ell > k$ such that $b_{\ell} = i$. The subarray 
$\begin{bmatrix}
    a_{j}  & a_{k} & a_{\ell} \\
    b_{j}  & b_{k} & b_{\ell}
\end{bmatrix} = \begin{bmatrix}
    i+1  & i+1 & a_{\ell} \\
    b_{j}  & b_{k} & i
\end{bmatrix}$
is then a valley as $b_{k} \leqslant i = b_{\ell} < i+1 = a_{k}$ and $b_{k} < b_{j}$. This contradicts $\mathcal{A}_G$ being a $PV$-free Burge array.

\smallskip

\noindent
\textbf{Case 2:} Assume that $a_{j} = i+1$ and $b_{k} = i+1$.
This once again implies that $b_{k}$ is the second leftmost $i+1$ in $\tilde{R}_{i}(\mathcal{A}_G)$. As $\mathcal{A}_G$ contains at least two $i$-pairs, 
there must exist $\ell > k$ such that $b_{\ell} = i$. Let $j < s < \ell$ be the leftmost index such that $b_{j} \leqslant i =b_{\ell} < b_{s}$. Note that such an 
$s$ exists as $k$ satisfies the desired conditions. The subarray 
$\begin{bmatrix}
    a_{j}  & a_{s} & a_{\ell} \\
    b_{j}  & b_{s} & b_{\ell}
\end{bmatrix} = \begin{bmatrix}
    i+1  & a_{s} & a_{\ell} \\
    b_{j}  & b_{s} & i
\end{bmatrix}$
is a peak as $b_{j} \leqslant i =b_{\ell} < b_{s}$ and $a_{j} = i+1 \leqslant b_{s}$. This contradicts $\mathcal{A}_G$ being a $PV$-free Burge array.

\smallskip

\noindent
\textbf{Case 3:} Assume that $b_{j} = i+1$ and $b_{k} = i+1$.
This implies that $b_{j}$ is the leftmost $i+1$ in $\tilde{R}_{i}(\mathcal{A}_G)$. To have at least two $i$-pairs, $\mathcal{A}_G$ must contain columns 
$\ell$ and $m$ such that $j< \ell < m$ and $b_{\ell} = b_{m} = i$. The subarray  $\begin{bmatrix}
    a_{j}  & a_{\ell} & a_{m} \\
    b_{j}  & b_{\ell} & b_{m}
\end{bmatrix} = \begin{bmatrix}
    a_{j}  & a_{\ell} & a_{m} \\
    i+1  & i & i
\end{bmatrix}$ 
is a valley as $b_{\ell} = i \leqslant i = b_{m} < a_{\ell}$ and $b_{\ell} = i < i+1 = b_{j}$. This contradicts $\mathcal{A}_G$ being a $PV$-free Burge array.
\end{proof}

\begin{lemma}\label{lemma.knuth_equivalence}
Let $T_{G}$ be the threshold tableau associated to $\mathcal{A}_G$ under the Burge correspondence. Then $\tilde{R}_{i}(\mathcal{A}_G)$ is Knuth 
equivalent to $R_{i}(T_G)$. 
\end{lemma}

\begin{proof}
Since $T_G$ is a hook tableau, $R_{i}(T_G)$ has at most one $i$-pair. Thus by Lemma~\ref{lemma.i_pair} it suffices to prove that 
$\tilde{R}_{i}(\mathcal{A}_G)$ has an $i$-pair if and only if 
$R_{i}(T_G)$ has an $i$-pair, since the content of $T_G$ and $\AG$ are the same.

Assume that $\tilde{R}_{i}(\mathcal{A}_G)$ has an $i$-pair. Let $k$ be the index of the column containing the $i+1$ that is in the $i$-pair. First 
assume that $a_{k} = i+1$. If $k = 1$, then $a_{k}$ is recorded in the leg of $T_G$ and $R_{i}(T_{G})$ has an $i$-pair. If $b_{k-1} \leqslant b_{k}$ 
then $b_{k}$, when inserted, does not bump an element. Otherwise it would bump an element greater than the element bumped by $b_{k-1}$ which 
would create a non-hook shape. Thus, $a_{k}$ is recorded in the leg of $T_G$ and $R_{i}(T_{G})$ has an $i$-pair. If $b_{k-1} > b_{k}$, then in order 
for $\tilde{R}_{i}(\mathcal{A}_G)$ to have an $i$-pair, there must exist $\ell > k$ such that $b_{\ell} = i$. The subarray 
$\begin{bmatrix}
    a_{k-1}  & a_{k} & a_{\ell} \\
    b_{k-1}  & b_{k} & b_{\ell}
\end{bmatrix} = \begin{bmatrix}
    a_{k-1}  & i+1 & a_{\ell} \\
    b_{k-1}  & b_{k} & i
\end{bmatrix}$ 
is a valley as $b_{k} \leqslant i = b_{\ell} < i+1 = a_{k}$ and $b_{k-1} > b_{k}$ which is a contradiction. Next assume that $b_{k} = i+1$. In order for 
$\tilde{R}_{i}(\mathcal{A}_G)$ to have an $i$-pair, there must exist $\ell > k$ such that $b_{\ell} = i$. This implies that some $i+1$ must be bumped 
into the leg by an insertion letter whose index is at most $\ell$. Thus, $T_G$ has an $i+1$ in its leg and $R_{i}(T_{G})$ has an $i$-pair. Thus, 
$R_{i}(T_G)$ has an $i$-pair whenever $\tilde{R}_{i}(\mathcal{A}_G)$ has an $i$-pair.

Assume that $R_{i}(T_G)$ has an $i$-pair. This implies that $T_{G}$ has an $i+1$ in its leg. If the $i+1$ in the leg corresponds to a recording letter 
$a_{j}$ for some $j$, then $b_{j}$ does not bump an element when inserted. Otherwise $a_{j} = i+1$ would get placed in the first row of $T_{G}$ and 
cannot be bumped into the leg as $a_{1}\leqslant a_{j}$. This implies that either $j=1$ or $b_{j-1} \leqslant b_{j}$. In either case, 
$\tilde{R}_{i}(\mathcal{A}_G)$ contains an $i$-pair. Assume the $i+1$ in the leg corresponds to an insertion letter $b_{j}$. Since $b_{j}$ must be 
bumped into the leg and be $i$-paired with some $i$, there exists $j< k$ such that $b_{k} = i$. This implies that $\tilde{R}_{i}(\mathcal{A}_G)$ contains 
an $i$-pair. Thus, $\tilde{R}_{i}(\mathcal{A}_G)$ has an $i$-pair whenever $R_{i}(T_G)$ has an $i$-pair.
\end{proof}

\begin{proposition}
\label{prop.hook_intertwine}
Let $\mathcal{A}_G$ be a PV-free Burge array and let $T_{G}$ be its associated threshold tableau.
\begin{enumerate}
\item
If $\tilde{f}_i(\mathcal{A}_{G}) \not = 0$, then $\tilde{f}_{i}(\mathcal{A}_{G}) = \mathcal{A}'_{G}$, where $\mathcal{A}'_{G}$ is the associated 
Burge array of $f_{i}(T_{G})$.
\item
If $\tilde{e}_i(\mathcal{A}_{G}) \not = 0$, then $\tilde{e}_{i}(\mathcal{A}_{G}) = \tilde{\mathcal{A}}_{G}$, where $\tilde{\mathcal{A}}_{G}$ is the 
associated Burge array of $e_{i}(T_{G})$.
\end{enumerate}
\end{proposition}

\begin{proof}
As $\tilde{f}_i$ and $\tilde{e}_{i}$ are clearly partial inverses, it suffices to just prove part (1). From Lemma~\ref{lemma.knuth_equivalence}, we 
have that $f_{i}(T_{G})$ is not annihilated. Let $s$ be the column index of the rightmost $i$ in $\AG$ and denote this rightmost $i$ 
by $\bar{i}$. We claim that $\bar{i}$ corresponds to the rightmost $i$ in $R_{i}(T_{G})$. If $b_{s} = \bar{i}$, then $\bar{i}$ is inserted into the arm 
to the right of all preexisting $i$'s when column $s$ is inserted and will remain the rightmost $i$ by the properties of Schensted row insertions. 
If $a_{s} = \bar{i}$ and $a_{s-1} = i$, then when column $s$ is inserted $b_{s}$ will bump an element from the arm. As $T_{G}$ is hook-shaped, 
$a_{s} = \bar{i}$ will be recorded into the arm to the right of all preexisting $i$'s and will remain the rightmost $i$. If $a_{s} = \bar{i}$ and 
$a_{s-1} \not = i$, then $a_{s}$ is the only $i$ in $\mathcal{A}_{G}$ and is trivially the rightmost $i$ in $R_{i}(T_{G})$ which proves our claim. 
As $f_{i}(T_{G})$ is not annihilated and $T_{G}$ is hook-shaped, $\bar{i}$ is the rightmost unpaired $i$ of $T_{G}$ and is changed to an $i+1$ by $f_{i}$.

Let $r$ be the column index of the last column of $\mathcal{A}_{G}$. We denote by $S_{t}$ the tableau obtained by reverse inserting columns 
$t+1$ through $r$ of $f_{i}(T_{G})$ and $T_{t}$ the tableau obtained by inserting columns $1$ through $t$ of $\mathcal{A}_{G}$. We first assume that $r>s$ 
and prove that columns $s+1$ through $r$ are the same in both arrays so that we may take $r$ to be the same as $s$ later in the proof. More specifically 
using induction, we will prove that for $s+1 \leqslant t \leqslant r$ column $t$ of $\mathcal{A}'_{G}$ is the same as column $t$ in $\mathcal{A}_{G}$ and 
$S_{t-1} = f_{i}(T_{t-1})$. 

Let $\ell$ and $a$ be the largest entries in the leg and arm of $T_{r} = T_{G}$, respectively. 
Note that $\ell$ and $a$ are also the largest entries in the leg and arm of $S_{r} = f_{i}(T_{G})$, respectively, as $r > s$. We break into subcases depending on whether $\ell > a$ or $\ell \leqslant a$. 

When $\ell > a$, the column $[\ell, a]^{T}$ is obtained by reverse inserting both $T_{r}$ and $S_{r}$ implying the $r$-th column of $\mathcal{A}'_{G}$ 
is the same as the $r$-th column of $\mathcal{A}_{G}$. We claim that $\ell > i+1$. This is clearly true if $\bar{i}$ is in the leg of $T_{G}$ as we assume $r > s$. 
If $\bar{i}$ is in the arm of $T_{r}$, then $a \geqslant i+1$ as $a$ is to the right of $\bar{i}$ and $\bar{i}$ is the rightmost $i$ in $T_{r}$. Hence the claim is true,
and $\bar{i}$ remains the rightmost unpaired $i$ in $T_{r-1}$. Thus, $f_{i}$ acts on $T_{r-1}$ by changing $\bar{i}$ to an $i+1$ which is precisely $S_{r-1}$. 

When $\ell \leqslant a$, $a$ is removed from $T_{r}$ and a number which we denote by $b_{r}$ is reverse inserted. This implies that the $r$-th column of
$\mathcal{A}_{G}$ is  $[a, b_{r}]^{T}$. Similarly, reverse inserting a column from $S_{r}$, $a$ is removed from  $S_{r}$ and a number which we denote by 
$c_{r}$ is reverse inserted. Note that $b_{r}$ and $c_{r}$ are equal and are in the same cell of their corresponding tableau except if an $i$ and $i+1$ lie 
in the arm of $S_{r}$ and the $(2,1)$ cell of $S_{r}$ is an $i+1$. If an $i$ and $i+1$ lie in the arm of $S_{r}$, then the value of the $(2,1)$ entry in 
$T_r$ cannot be $\bar{i}$ as this would contradict $\bar{i}$ being the rightmost unpaired $i$ of $T_{r}$.  Thus, if an $i$ and $i+1$ lie in the 
arm of $S_{r}$ and the $(2,1)$ cell of $S_{r}$ is an $i+1$, then the $(2,1)$ cell of $T_{r}$ is an $i+1$. However, this would imply $b_{r} = i$ 
as it would need to bump an $i+1$ so that $a$ is in the arm of $T_{r}$, an $i+1$ is in the $(2,1)$ cell, and the $i$ in the arm is not bumped. 
This contradicts $r > s$. Therefore the $r$-th column of $\mathcal{A}_{G}$ and $\mathcal{A}'_{G}$ match. As $b_{r} \not = i$, $\bar{i}$ is still 
the rightmost unpaired $i$ of $T_{r-1}$ and $f_{i}(T_{r-1}) = S_{r-1}$. Assume that column $t$ of $\mathcal{A}'_{G}$ is the same as column 
$t$ in $\mathcal{A}_{G}$ and $S_{t-1} = f_{i}(T_{t-1})$ for some $s+1 < t \leqslant r$. Repeating the argument in the base case, we see that this 
holds for the case $t-1$ as well.

From the definition of $\tilde{f}_{i}$, we have that the columns $s+1$ through $r$ of $\tilde{f}_{i}(\mathcal{A}_{G})$ are the same as the columns 
$s+1$ through $r$ of $\mathcal{A}_{G}$. By the previous paragraph, these columns are also equal to columns $s+1$ through $r$ of 
$\mathcal{A}'_{G}$. We now assume that $r = s$ and prove that the columns $1$ through $s$ of $\tilde{f}_{i}(\mathcal{A}_{G})$ and $\mathcal{A}'_{G}$ 
are equal by breaking into cases based off the position of $\bar{i}$ in $T_{s}$.

\smallskip

\noindent
\textbf{Case 1:} $\bar{i}$ is in the leg of $T_{s}$.

Since $\bar{i}$ is the rightmost unpaired $i$ of $T_{s}$, this implies that $T_{s}$ does 
not contain any other $i$ except for $\bar{i}$.  Also in order for $\bar{i}$ to be in column $s$ of $\mathcal{A}_{G}$, we must have $a_{s} = \bar{i}$. 
Thus, the column obtained from $T_{s}$ by reverse inserting is of the form $[\bar{i}, a]^{T}$ where $a$ is the largest value in the arm of $T_{s}$. 
In $S_{s}$, $\bar{i}$ is replaced with an $i+1$. Thus, the column obtained from $S_{s}$ by reverse inserting is of the form $[i+1, a]^{T}$. We see 
that $T_{s-1}$ and $S_{s-1}$ are equal implying columns $1$ through $s-1$ of $\mathcal{A}_{G}$ and $\mathcal{A}'_{G}$ are equal. Note that 
as $\bar{i}$ is the only $i$ in $\mathcal{A}_{G}$ and $a_{s} = \bar{i}$, we have $\tilde{f}_{i}$ acts by changing $\bar{i}$ to an $i+1$ in
 $\mathcal{A}_{G}$. This is precisely the form of $\mathcal{A}'_{G}$ implying $\tilde{f}_{i}(\mathcal{A}_{G}) = \mathcal{A}'_{G}$.

\smallskip

\noindent
\textbf{Case 2:} $\bar{i}$ is in the arm of $T_{s}$.

Let $\ell$ and $a$ be the largest elements in the leg and arm of $T_{s}$, respectively. We break into subcases.

Assume $\bar{i} \geqslant \ell$. In order for $\bar{i}$ to be in the $s$-th column of $\mathcal{A}_{G}$, $\bar{i}$ must be $a$; otherwise the column 
reverse inserted from $T_{s}$ would be $[a, \bar{i}]^{T}$. This implies $\bar{i}$ must bump an entry in the arm of $T_{s-1}$ into its leg which 
would contradict $\bar{i} \geqslant \ell$.  In particular, this implies that $a= \bar{i} \geqslant \ell$.
Furthermore, the column obtained from $T_{s}$ when reverse inserting is of the form $[\bar{i}, b_{s}]^{T}$, 
where $b_{s}$ is the largest element in the arm of $T_{s}$ strictly less than the entry in the cell $(2,1)$. Since $\bar{i}$ is the largest value in 
$T_{s}$, the $i+1$ in $S_{s}$ created by applying $f_{i}$ to $T_{s}$ is also the largest value. Thus, the column obtained from $S_{s}$ when reverse inserting is of the form 
$[i+1, c_{s}]^{T}$, where $c_{s}$ is the largest element in the arm of $S_{s}$ strictly less than the entry in the cell $(2,1)$. Since $f_{i}(T_{s}) = S_{s}$, 
we have $b_{s} = c_{s}$ and $T_{s-1} = S_{s-1}$. Thus, columns $1$ through $s-1$ of $\mathcal{A}_{G}$ are identical to the corresponding columns 
of $\mathcal{A}'_{G}$. As $\bar{i}$ is the rightmost $i$ in $\mathcal{A}_{G}$ and is in the top row, it is also the rightmost unpaired $i$ in $\mathcal{A}_{G}$. 
Thus, $\tilde{f}_{i}$ acts by changing $\bar{i}$ to an $i+1$ in $\mathcal{A}_{G}$. Thus, $\tilde{f}_{i}(\mathcal{A}_{G}) = \mathcal{A}'_{G}$.

Assume now that $\ell > a$ and $\ell \not = i+1$. In order for $\bar{i}$ to be in the $s$-th column of $\mathcal{A}_{G}$, $\bar{i}$ must be equal to $a$. When reversing the Burge correspondence of $T_{s}$, the column obtained is then of the form $[\ell, \bar{i}]^{T}$. 
From our assumption, we also have $\ell > i+1$ which implies that the column obtained from $S_{s}$ when reversing the Burge correspondence is 
$[\ell, i+1]^{T}$. Observe that $T_{s-1} = S_{s-1}$ which forces columns $1$ through $s-1$ of $\mathcal{A}_{G}$ and $\mathcal{A}'_{G}$ to be equal. 
We prove that $\bar{i}$ is the rightmost unpaired $i$ in $\mathcal{A}_{G}$. For there to be a hope that this claim is not true, then there must exist either a column $[i+1, b_{j}]^{T}$ with $b_{j} \not = i$ or 
$[a_{j}, i+1]^{T}$ in $\mathcal{A}_{G}$. Note that there can only be one column with an $i+1$ in the top row; otherwise it would form a valley with columns 
$[i+1, b_{j}]^{T}$ and $[\ell, \bar{i}]^{T}$. If there exists a column of the form $[i+1, b_{j}]^{T}$, then in order for $\bar{i}$ to be unpaired in $T_{s}$ there 
must be an $i$ in a column $j+1$ through $s-1$ or a column of the form $[i, b_{j'}]^{T}$. Note that if there is an $i$ in columns $j+1$ through $s-1$ this would imply $\bar{i}$ is the rightmost unpaired $i$ of $\AG$.  If there is no $i$ in columns $j+1$ through $s-1$,
we must then have that column $j-1$ is of the form $[i, b_{j-1}]^{T}$. If $b_{j-1} > b_{j}$, then 
 $\begin{bmatrix}
    i  & i+1 & \ell \\
    b_{j-1}  & b_{j} & \bar{i}
\end{bmatrix}$ 
is a valley implying $b_{j-1} \leqslant b_{j}$. Thus, if there exists a column of the form $[i+1, b_{j}]^{T}$ with no $i$'s in columns $j$ through $s-1$, 
then $\bar{i}$ must be the rightmost unpaired $i$. If there exists a column of the form $[a_{j}, i+1]^{T}$ with no $i$'s in columns $j$ through $s-1$, 
this would imply $\bar{i}$ when inserted would bump an element from the arm of $T_{s}$ contradicting that $\ell$ is in the leg of $T_{s}$. Thus, if there exists a column of the form $[a_{j}, i+1]^{T}$, then there exists an $i$ in columns $j$ through $s-1$. Therefore, 
$\bar{i}$ is the rightmost unpaired $i$ of $\mathcal{A}_{G}$ and $\tilde{f}_{i}$ acts by changing $\bar{i}$ to an $i+1$ in the $s$-th column. This 
is precisely $\mathcal{A}'_{G}$.

Assume that $\ell > a$ and $\ell = i+1$. In order for $\bar{i}$ to be in column $s$ of $\mathcal{A}_{G}$, $a$ must be precisely $\bar{i}$. When 
reversing the Burge correspondence of $T_{s}$, the column obtained is then of the form $[\ell = i+1, \bar{i}]^{T}$. Since $\ell = i+1$ came from 
the leg of $T_{s}$, there must exist an $i$ somewhere in columns $1$ through $s-1$; otherwise $\bar{i}$ would not be the rightmost 
unpaired $i$ in $T_{s}$. Thus, $\bar{i}$ is also the rightmost unpaired $i$ of $\mathcal{A}_{G}$.  Let $x$ be the value in the $(2,1)$ position 
of $T_{s}$ and let $y$ be the rightmost element in the arm of $T_{s}$ that is strictly less than $x$. Note that $x < i+1 = \ell$; otherwise 
$T_{s}$ would have shape $(1,1)$ and $R_{i}(T_{s}) = i+1 \, i$. This implies $y \not = i$. We also have $y$ is not an element 
in the top row of $\AG$. Otherwise if $a_{m} = y$ for some $m$, then $x > y \geqslant a_{1}$ which contradicts $x$ being in the $(2,1)$ cell. 
Thus, $y = b_{m}$ for some $1 \leqslant m \leqslant s-1$. 

Assume now that $x$ is a recording letter. As $x$ lies in the $(2,1)$ cell of $T_{s}$, we must 
have $x = a_{1}$. This implies $b_{1} \leqslant b_{2} \leqslant \cdots \leqslant b_{s} = \bar{i}$ and $a_{1} < a_{2} < \cdots < a_{s} = \ell$. Recall that in 
$S_{s}$, $b_{s} = \bar{i}$ is replaced with an $i+1$. Hence, when reverse inserting a column from $S_{s}$, the $i+1$ that replaced $\bar{i}$ is removed, 
$y = b_{m}$ is replaced by 
$x = a_{1}$, and the rest of the entries in the leg are shifted up. Thus column $s$ of $\mathcal{A}'_{G}$ is of the form $[i+1, b_{m}]^{T}$. 
We see that $S_{s-1}$ is then associated to the Burge array 
$\begin{bmatrix}
    a_{2}  & a_{3} & \dots & a_{m} & a_{m+1} & a_{m+2} & \ldots & a_{s} = \ell = i+1 \\
    b_{1}  & b_{2} & \dots & b_{m-1} & a_{1} & b_{m+1} & \ldots & b_{s-1}
\end{bmatrix}$ 
which mimics the action of $\tilde{f}_{i}$ on $\AG$ as $b_{1} \leqslant \cdots \leqslant b_{s-1}$ and $b_{m}$ is the rightmost entry such that 
$b_{m} < a_{1}$. 

Assume now that $x$ lies in the bottom row of $\mathcal{A}_{G}$. This implies that there exists $b_{n}$ that bumped $x$ in 
$T_{n-1}$ to the $(2,1)$ cell. Note that $b_{n-1} > b_{n}$; otherwise 
$\begin{bmatrix}
    a_{z}  & a_{n-1} & a_{n} \\
    x  & b_{n-1} & b_{n}
\end{bmatrix}$ 
would form a valley in $\mathcal{A}_{G}$. Since $x$ is the value in the $(2,1)$ cell, $b_{n} \leqslant b_{n+1} \leqslant \cdots \leqslant b_{s} = \bar{i}$ in 
$\mathcal{A}_{G}$. We also have $b_{n+1} \geqslant a_{n}$. Otherwise 
$\begin{bmatrix}
    a_{z}  & a_{n} & a_{n+1} \\
    x  & b_{n} & b_{n+1}
\end{bmatrix}$ 
would form a valley. Thus, $b_{n}$ is the largest element in the bottom row from column $n$ to $s-1$ such that $b_{n} < a_{n}$ which implies 
$m = n$. Moreover, $a_{n} < \cdots < a_{s} = \ell$ are all in the legs of both of $T_{s}$ and $S_{s}$. As $S_{s}$ differs from $T_{s}$ by changing 
$\bar{i}$ to an $i+1$, we have that reversing the Burge correspondence removes the $i+1$ from the arm of $S_{s}$, $b_{n}$ is replaced by $x$, 
and the rest of the leg entries are shifted up. Thus column $s$ of $\mathcal{A}'_{G}$ is $[i+1, b_{n}]^{T}$. As $b_{n}$ bumped $x$ out, we see that 
$x$ is in the cell that it was originally inserted into. We see that reverse the Burge correspondence for $S_{s-1}$ up to $S_{n-1}$, we 
get the columns 
$\begin{bmatrix}
    a_{n+1}  & a_{n+2} & \dots & a_{s-1} & a_{s} = \ell = i+1 \\
    a_{n}  & b_{n+1} & \dots & b_{s-2} & b_{s-1}
\end{bmatrix}$
and $S_{n-1}$ is equal to $T_{n-1}$ as $x$ is in its original cell. These changes to $\mathcal{A}_{G}$ are seen to be the same as $\tilde{f}_{i}$ 
as $n$ is the leftmost column index such that $b_{n} \leqslant \cdots \leqslant b_{s-1}$ and $n$ is the rightmost column index between $n$ and 
$s-1$ such that $b_{n} < a_{n}$.

Assume that $\bar{i} < \ell \leqslant a$ and the $(2,1)$ cell of $T_{s}$ is not an $i+1$. Let $x$ be the value in the $(2,1)$ position of $T_{s}$ 
and let $y$ be the rightmost element in the arm of $T_{s}$ that is strictly less than $x$. Note that for $\bar{i}$ to be in the $s$-th column of 
$\mathcal{A}_{G}$, $y$ must equal $\bar{i}$. As $x$ is not equal to $i+1$, we have also $i+1 < x$. When reversing the Burge correspondence 
of $T_{s}$ and $S_{s}$ the columns obtained are $[a, \bar{i}]^{T}$ and $[a, i+1]^{T}$ respectively where the $i+1$ reverse bumped from 
$S_{s}$ was the $i+1$ created by $f_{i}$. We have $T_{s-1}$ and $S_{s-1}$ are equal implying columns $1$ through $s-1$ of 
$\mathcal{A}_{G}$ and $\mathcal{A}'_{G}$ are the same. Note that there cannot be an $i+1$ in columns $1$ through $s-1$. Otherwise 
there would either be an $i+1$ in the leg of $T_{s-1}$ which would contradict $x > i+1$ or $i+1$ is in the arm of $T_{s-1}$ in which case 
$\bar{i}$ would bump an $i+1$ instead of $x$. Thus, $\bar{i}$ is the rightmost unpaired $i$ of $\mathcal{A}_{G}$ and $\tilde{f}_{i}$ changes 
$\bar{i}$ to an $i+1$.

Assume that $\bar{i} < \ell \leqslant a$ and the $(2,1)$ cell of $T_{s}$ is an $i+1$. Let $x = i+1$ be the value in the $(2,1)$ cell of $T_{s}$ 
and $y$ be the rightmost element in the arm of $T_{s}$ that is strictly less than $x$. For $\bar{i}$ to be in the $s$-th column of 
$\mathcal{A}_{G}$, $y$ must equal $\bar{i}$. As column $s$ in $\mathcal{A}_{G}$ is of the form $[a, \bar{i}]^{T}$, $\bar{i}$ bumps 
$i+1$ from the arm of $T_{s-1}$. Since $x$ is bumped into the cell $(2,1)$, we have $x = b_{n}$ for some $n$. Observe that no $i$ 
can be strictly between columns $n$ through $s$ in $\AG$; otherwise $T_{s-1}$ would contain an $i+1$ in its leg. From this observation and the fact that 
$\bar{i}$ is the rightmost unpaired $i$ in $T_{G}$, there must exist an $i$ somewhere in columns $1$ through $n-1$ in 
$\mathcal{A}_{G}$. 

Let $m$ be the column of the index of the second rightmost $i$ in $\mathcal{A}_{G}$ which by the reasoning 
above satisfies $m < n$. We claim that $i$ must be the insertion letter in column $m$, i.e. $b_m = i$. If $a_{m} = i$, then $b_{m}$ must have bumped 
an element in $T_{m-1}$ when inserted; otherwise an $i$ would be present in the leg of $T_{s-1}$. Let $b_{z}$ be the element bumped by $b_{m}$. 
This implies $b_{z}$ 
is in the leg of $T_{m-1}$; however, $b_{z} < a_{z} \leqslant a_{m} = i < i+1$ which would be a contradiction. Thus our claim that $b_{m} = i$ 
holds. We also have $a_{m} \not = i+1$ as $\tilde{R}_{i}(\mathcal{A}_{G})$ can have at most one $i$-pair by Lemma~\ref{lemma.i_pair}. 
From these two facts, $b_{m}$ can not have bumped an element $b_{u}$ when inserted into $T_{m-1}$. Otherwise 
$\begin{bmatrix}
    a_{u}  & a_{m}  & a \\
    b_{u}  & b_{m}  & \bar{i}
\end{bmatrix}$ 
would be a valley. As $y = \bar{i}$ is turned into an $i+1$ in $S_{s}$, the rightmost element in the arm of $S_{s}$ that is strictly less 
than $x = i+1$ is $b_{m}$. Thus, the column reverse inserted from $S_{s}$ is $[a, b_{m} = i]^{T}$ while the column reverse inserted from 
$T_{s}$ is $[a, \bar{i}]^{T}$. This implies the $s$-th columns of $\mathcal{A}_{G}$ and $\mathcal{A}'_{G}$ are the same, but $S_{s-1}$ 
differs from $T_{s-1}$. Note that columns $m+1$ through $s-1$ of both $\mathcal{A}_{G}$ and $\mathcal{A}'_{G}$ are the same as the 
reverse insertions of $T_{s}$ and $S_{s}$ in these steps do not involve $b_{m}$ or the $i+1$ created by $f_{i}$ respectively. As $b_{m}$ 
did not bump an element when inserted into $T_{m-1}$, we have $a_{m}$ is in the leg of both $T_{m}$ and $S_{m}$ and is the largest entry. 
We see $[a_{m}, b_{m}]^{T}$ is reverse inserted from $T_{m}$ and $[a_{m}, i+1]^{T}$ is reverse inserted from $S_{m}$ and $S_{s-1} = T_{s-1}$. 
Thus, the only difference between $\mathcal{A}_{G}$ and $\mathcal{A}'_{G}$ is that $b_{m}$ is changed to an $i+1$ in $\mathcal{A}_{G}$. 
Note that $b_{m}$ is the rightmost unpaired $i$ in $\mathcal{A}_{G}$ since it is the second rightmost $i$ and $\bar{i}$ is $i$-paired with 
the $b_{n}$. Thus, $\tilde{f}_{i}$ acts by changing $b_{m}$ in $\mathcal{A}_{G}$ to an $i+1$ which is precisely $\mathcal{A}'_{G}$.
\end{proof}

\begin{corollary}
Let $C_{m}$ be the set of all $PV$-free Burge arrays with entries at most $m$. Then $C_{m}$ together with the operators $\tilde{f}_{i}$ and $\tilde{e}_{i}$ 
forms a Stembridge crystal of type $A_{m-1}$.
\end{corollary}

\begin{proof}
The Burge correspondence is a crystal isomorphism between $C_{m}$ and 
$\displaystyle \bigsqcup_{\substack{\lambda \text{ hook-shaped,} \\ \text{threshold}}} \Tab_m(\lambda)$ by 
Proposition~\ref{prop.hook_intertwine}, where $\Tab_m(\lambda)$ is the set of all semistandard Young tableaux of shape $\lambda$ 
and entries at most $m$ together with the usual crystal operators as in Definition~\ref{definition.crystal}. Since $\Tab_m(\lambda)$
forms a Stembridge crystal, so does $C_m$.
\end{proof}

\begin{corollary} 
Let $\mathcal{A}_{G}$ be a PV-free Burge array corresponding to a graph $G$ on $n$ vertices. Then $\mathcal{A}_{G}$ is highest 
weight if and only if $G$ is star-shaped (up to singletons) such that the central vertex is labelled $1$ and the other vertices have labels 
$\{2, \ldots, n\}$.
\end{corollary}

For a crystal $C$, an element $b\in C$ is called \defn{extremal} if either $f_i(b) = 0$ or $e_i(b)=0$ for each $i$ in the index set and its weight
is in the Weyl orbit of the highest weight element in the crystal component (see~\cite{Kashiwara.1994}). 
Let the weight of the highest weight vector $u\in C$ (which satisfies $e_i(u)=0$ for all $i$) be $\wt(u)=\lambda$. The weight of the extremal vectors are 
permutations of $\lambda$. The tableaux under the Burge correspondence are threshold shapes. Hence, by the definition of threshold graphs,
the extremal vectors of the crystal  correspond to threshold graphs under the Burge correspondence.

%%%%%%%%%%%%%%%%%%%%%%%%%%%%%%%%%%%%%
\bibliographystyle{plain}
\bibliography{erdos_gallai}{}

\end{document}